\documentclass[12pt]{article}
\usepackage{amsmath}
\usepackage{amssymb}
\allowdisplaybreaks

\def\q{\quad}
\def\qq{\qquad}

\def\mod#1{\ (\text{\rm mod}\ #1)}
\def\t{\text}
\def\f{\frac}
\def\e{\equiv}
\def\b{\binom}
\def\sls#1#2{(\f{#1}{#2})}
 \def\ls#1#2{\big(\f{#1}{#2}\big)}
\def\Ls#1#2{\Big(\f{#1}{#2}\Big)}
\def\ap{\langle a\rangle_p}
\def\xp{\langle x\rangle_p}

\begin{document}
 \centerline {\bf
Generalizations of the Christoffel-Darboux formula}
\centerline {\bf and congruences involving Ap\'ery-like numbers}
\par\q\newline
\centerline{Zhi-Hong Sun}\newline \centerline{School of Mathematics
and Statistics} \centerline{Huaiyin Normal University}
\centerline{Huaian, Jiangsu 223300, P.R. China} \centerline{Email:
zhsun@hytc.edu.cn} \centerline{URL:
https://maths.hytc.edu.cn/szh.htm}
\vskip0.2cm
\par\q
\par {\bf Abstract.} In this paper, we first extend the  Christoffel-Darboux formula for orthogonal polynomials
to general three-term recurrence sequences, and then investigate the identities and congruences for $g_n(x)$ and $v_n(x)$ given by
\begin{align*}
&g_0(x)=1,\ g_1(x)=\frac{x+1}2,\ (n+1)^2g_{n+1}(x)=\Big(2n(n+1)+\frac{x+1}2\Big)g_n(x)-n^2g_{n-1}(x)\ (n\ge 1),
\\&v_0(x)=1,\ v_1(x)=x,\ (n+1)^3v_{n+1}(x)=(2n+1)(n(n+1)+x)v_n(x)-n^3v_{n-1}(x)\ (n\ge 1).\end{align*}
\par\q

\section*{1. Introduction}
\par Let $G_n(x)$ and $V_n(x)$ be defined by
\begin{align*}
&G_n(x)=\sum_{k=0}^n\b nk(-1)^k\b xk\b{-1-x}k,
\\&V_n(x)=\sum_{k=0}^n\b nk\b{n+k}k(-1)^k\b
xk\b{-1-x}k\ (n=0,1,2,\ldots).\end{align*}
Then
\begin{align*}
&G_0(x)=1,\q G_1(x)=x^2+x+1,\ \\&(n+1)^2G_{n+1}(x)=(2n(n+1)+x^2+x+1)G_n(x)-n^2G_{n-1}(x)\ (n\ge 1),
\\&V_0(x)=1,\q V_1(x)=2x^2+2x+1,
\\&(n+1)^3V_{n+1}(x)=(2n+1)(n(n+1)+2x^2+2x+1)V_n(x)-n^3V_{n-1}(x)\ (n\ge 1).\end{align*}
by [15, (22.10)] and [14]. Thus,
\begin{align*}
&(n+1)^2\cdot m^{n+1}G_{n+1}(x)
\\&=(2mn(n+1)+m(x^2+x+1))\cdot m^nG_n(x)
-m^2n^2\cdot m^{n-1}G_{n-1}(x)\ (n\ge 1),
\\&(n+1)^3\cdot m^{n+1}V_{n+1}(x)
\\&=(2n+1)(mn(n+1)+m(2x^2+2x+1))\cdot m^nV_n(x)
-m^2n^3\cdot m^{n-1}V_{n-1}(x)\ (n\ge 1).
\end{align*}
From [14,15],
\begin{align*}
&G_n(x)=\sum_{k=0}^n\b xk^2(-1)^{n-k}\b{-1-x}{n-k},
\\&V_n(x)=\sum_{k=0}^n\b nk\b{n+k}k(-1)^{n-k}G_k(x)
=\sum_{k=0}^n\b xk^2\b{-1-x}{n-k}^2.\tag{1.1}
\end{align*}
By [16],
$$\sum_{n=0}^{p-1}\f{V_n(x)}{m^n}\e
\Big(\sum_{k=0}^{p-1}\f{G_k(x)}{m^k}\Big)^2\mod p\eqno{(1.2)}$$ for any
prime $p>3$ and $m,x\in\Bbb Z_p$ with $m\not\e 0\mod p$, where $\Bbb
Z_p$ is the set of those rational numbers whose denominators are not divisible by $p$.
\par Let $V_n,V_n^{(3)},V_n^{(4)}$ and $V_n^{(6)}$ be defined by
\begin{align*} &V_n=16^nV_n\Big(-\f 12\Big)=\sum_{k=0}^n
\b nk\b{n+k}k(-1)^k\b{2k}k^216^{n-k},
\\&V_n^{(3)}=27^nV_n\Big(-\f 13\Big)=\sum_{k=0}^n\b
nk\b{n+k}k(-1)^k\b{2k}k\b{3k}k27^{n-k},
\\&V_n^{(4)}=64^nV_n\Big(-\f 14\Big)=\sum_{k=0}^n\b
nk\b{n+k}k(-1)^k\b{2k}k\b{4k}{2k}64^{n-k},
\\&V_n^{(6)}=432^nV_n\Big(-\f 16\Big)=\sum_{k=0}^n\b
nk\b{n+k}k(-1)^k\b{3k}k\b{6k}{3k}432^{n-k}.
\end{align*}
Then $V_n,V_n^{(3)},V_n^{(4)}$ and $V_n^{(6)}$ are Ap\'ery-like numbers of the first kind. Let $G_n,G_n^{(3)},G_n^{(4)}$ and $G_n^{(6)}$ be given by
\begin{align*} &G_n=16^nG_n\Big(-\f 12\Big)=\sum_{k=0}^n\b
nk(-1)^k\b{2k}k^216^{n-k},
\\&G_n^{(3)}=27^nG_n\Big(-\f 13\Big)=\sum_{k=0}^n\b
nk(-1)^k\b{2k}k\b{3k}k27^{n-k},
\\&G_n^{(4)}=64^nG_n\Big(-\f 14\Big)=\sum_{k=0}^n\b
nk(-1)^k\b{2k}k\b{4k}{2k}64^{n-k},
\\&G_n^{(6)}=432^nG_n\Big(-\f 16\Big)=\sum_{k=0}^n\b
nk(-1)^k\b{3k}k\b{6k}{3k}432^{n-k}.\end{align*}
Then
$G_n,G_n^{(3)},G_n^{(4)}$ and $G_n^{(6)}$ are Ap\'ery-like numbers of the second kind. For the discussions on two kinds of Ap\'ery-like numbers see [4,6,7,10,13,14,16,19].
\par
In this paper, we introduce $g_n(x)$ and  $v_n(x)$ as generalizations of $G_n(x)$ and $V_n(x)$. Let $g_n(x)$ and $v_n(x)$ be given by
\begin{align*}
&g_0(x)=1,\ g_1(x)=\f{x+1}2,\ (n+1)^2g_{n+1}(x)=\Big(2n(n+1)+\f{x+1}2\Big)g_n(x)-n^2g_{n-1}(x)\ (n\ge 1),
\\&v_0(x)=1,\ v_1(x)=x,\ (n+1)^3v_{n+1}(x)=(2n+1)(n(n+1)+x)v_n(x)-n^3v_{n-1}(x)\ (n\ge 1).\end{align*}
Then clearly $G_n(x)=g_n(2x^2+2x+1)$ and $V_n(x)=v_n(2x^2+2x+1)$.
\par Let $\{p_n(x)\}$ be the orthogonal polynomials with the weight function $w(x)$ over the interval $(a,b)$. Suppose that
$k_n$ is the coefficient of $x^n$ in $p_n(x)$ and
$h_n=\int_a^b w(x)p_n(x)^2dx$. The Christoffel-Darboux formula (see [1,2,5]) states that
$$\sum_{m=0}^n\f {p_m(x)p_m(y)}{h_m}=\f{k_n(p_{n+1}(x)p_n(y)
-p_n(x)p_{n+1}(y))}{h_nk_{n+1}(x-y)}.$$
Let $y\to x$. We get
$$\sum_{m=0}^n\f 1{h_m}p_m(x)^2
=\f{k_n}{h_nk_{n+1}}(p_n(x)p'_{n+1}(x)-p'_n(x)p_{n+1}(x)),$$
where $p'_n(x)=\f {d\,p_n(x)}{dx}$.
\par The paper is organized as follows. In Section 2, we extend the  Christoffel-Darboux formula. Suppose that $a_0=0$, $c\not=0$ and
$a_{n+1}u_{n+1}(x)=b_n(x+d_n)u_n(x)-ca_nu_{n-1}(x)\ (n=1,2,3,\ldots).$
We establish the closed formulas for
$$\sum_{n=0}^{p-1}\f{b_n}{c^n}u_n(x)u_n(y),\
\sum_{n=0}^{p-1}\f{b_n}{c^n}u_n(x)^2,\ \sum_{n=0}^{p-1}\f{b_n(x+y+2d_n)}{(-c)^n}u_n(x)u_n(y).$$
and
$$\sum_{n=0}^{p-1}\big(x-y+(-1)^{p-1-n}(x+y+2d_n)\big)\f{b_n}{c^n}
u_n(x)u_n(y).$$
Section 3 is devoted to giving the identities involving $g_n(x)$ and $v_n(x)$. As consequences of the identities in Section 2, we have the formulas for
$$\sum_{n=0}^{p-1}g_n(x)g_n(y)\q \sum_{n=0}^{p-1}g_n(x)^2,
\ \sum_{n=0}^{p-1}(2n+1)v_n(x)v_n(y)\ \t{and}\ \sum_{n=0}^{p-1}(2n+1)v_n(x)^2.$$
As applications, we deduce the closed formulas for
$\sum_{n=0}^{p-1}(2n+1)^{2r+1}v_n(x)$ and $\sum_{n=0}^{p-1}n^rg_n(x)$, where $p$ is a positive integer and $r\in\{0,1,2\}$. Actually the formulas for $r\ge 3$ can be deduced similarly.
\par In Section 4, using the identities in Section 3 we deduce the congruences for $\sum_{n=0}^{p-1}(2n+1)^mv_n(x)$ modulo $p^5$, where $p$ is an odd prime, $m\in\{1,3,5\}$, $x\in\Bbb Z_p$ and $\ls{2x-1}p=1$. Here $\sls ap$ is the Legendre symbol. Moreover, we deduce the congruences for $\sum_{n=0}^{p-1}(2n+1)^3\f{V_n}{16^n}$ and
$\sum_{n=0}^{p-1}(2n+1)^5\f{V_n}{16^n}$ modulo $p^7$. We also establish the congruence for $\sum_{n=0}^{p-1}(2n+1)v_n(x)^2$ modulo $p^4$. In particular, since $V_n(x)=v_n(2x^2+2x+1)$ we have
  $$\sum_{n=0}^{p-1}(2n+1)V_n(x)^2\e \f{1+2x'}{1+2x}p^2\mod
{p^4}$$
for $x\in\Bbb Z_p$ with $x\not\e -\f 12\mod p$,
where $x'=(x-\xp)/p$ and $\xp$ is given by
 $x\e \xp\mod p$ and $\xp\in\{0,1,\ldots,p-1\}$.
It should be mentioned that Yang and Liu[18] proved the congruence  $\sum_{n=0}^{p-1}(2n+1)V_n(x)^2\e 0\mod {p^2}$ for $x\not\e -\f 12\mod p$, which was conjectured by the author in [16]. We also note that the author
gave the congruence for $\sum_{n=0}^{p-1}V_n(x)$ modulo $p^4$ in [16], and Mao and Yang[8] obtained the congruences for $\sum_{n=0}^{p-1}(2n+1)^3V_n(x)$ and $\sum_{n=0}^{p-1}(2n+1)^5V_n(x)$ modulo $p^4$.
\par Let $p$ be an odd prime, $x\in\Bbb Z_p$ and $\sls{2x-1}p=1$. In Section 5, we establish the congruence
 for $\sum_{n=0}^{p-1}g_n(x)^2$ modulo $p^3$. As a consequence, we confirm the conjecture
  $$\sum_{n=0}^{p-1}G_n(x)^2\e (-1)^{\xp}\f{1+2x'}{1+2x}p\mod
{p^3}$$
due to the author's brother Z.W. Sun[17], where $x\in\Bbb Z_p$ satisfying $x\not\e -\f 12\mod p$.
   For $x\not\e 1,5,13\mod p$ and $r\in\{0,1,2\}$, the congruences for $\sum_{n=0}^{p-1}n^rg_n(x)$ modulo $p^4$ are also deduced. Furthermore, we establish the congruences for $\sum_{n=0}^{p-1}n^rG_n(-\f 1m)$ modulo $p^5$, where $m\in\{2,3,4,6\}$. Finally we pose the conjecture for $G_{\f{p-1}2}$ modulo $p^3$. We note that
the author gave the congruences for $\sum_{n=0}^{p-1}G_n(x)$ and $\sum_{n=0}^{p-1}nG_n(x)$ modulo $p^3$ in [14], and Mao and Yang [8] obtained the congruences for $\sum_{n=0}^{p-1}n^2G_n(x)$ and $\sum_{n=0}^{p-1}n^3G_n(x)$ modulo $p^3$.
\par In addition to the above notation, throughout this paper
we use the following notations: $\Bbb Z^+\f{\q}{\q}$the set of positive intgers, $\Bbb R\f{\q}{\q}$the set of real numbers, $[x]\f{\q}{\q}$the greatest integer not exceeding $x$, $H_0=0$, $H_n=\sum_{k=1}^n\f 1k\ (n\ge 1)$, $q_p(a)=(a^{p-1}-1)/p$. The Bernoulli numbers $\{B_n\}$, Euler numbers $\{E_n\}$ and the sequence
$\{U_n\}$ are defined by
\begin{align*} &B_0=1,\q\sum_{k=0}^{n-1}\b nkB_k=0\q(n\ge 2),\\& E_{2n-1}=0,\q E_0=1,\q E_{2n}=-\sum_{k=0}^{n-1}
\b {2n}{2k}E_{2k}\q(n\ge 1),
\\& U_{2n-1}=0,\q U_0=1,\q U_{2n}=-2\sum_{k=0}^{n-1}
\b {2n}{2k}U_{2k}\q(n\ge 1).
\end{align*}

\section*{2. Generalizations of the Christoffel-Darboux formula}
\par\q In this section, we generalize the Christoffel-Darboux formula for orthogonal polynomials.
\vskip0.2cm
\par\q{\bf Theorem 2.1} {\sl Suppose that $p\in\Bbb Z^+$, $a_0=0$, $c\not=0$ and
$$a_{n+1}u_{n+1}(x)=b_n(x+d_n)u_n(x)-ca_nu_{n-1}(x)\ (n=1,2,3,\ldots).$$
\par $(\t{\rm i})$ For $x\not=y$,
$$\sum_{n=0}^{p-1}\f{b_n}{c^n}u_n(x)u_n(y)=\f{a_p}{c^{p-1}(x-y)}
\big(u_p(x)u_{p-1}(y)-u_{p-1}(x)u_p(y)\big).$$
\par $(\t{\rm ii})$ We have}
$$\sum_{n=0}^{p-1}\f{b_n}{c^n}u_n(x)^2
=\f{a_p}{c^{p-1}}\big(u_{p-1}(x)u_p'(x)-u_p(x)u_{p-1}'(x)\big).$$

\par{\it Proof.} Set $h_n(x,y)=\f{a_n}{c^{n-1}}(u_n(x)u_{n-1}(y)-u_{n-1}(x)u_n(y))$.
Then
\begin{align*}h_{n+1}(x,y)&=\f{a_{n+1}}{c^{n}}u_{n+1}(x)u_n(y)-u_n(x)
\f{a_{n+1}}{c^{n}}u_{n+1}(y)
\\&=\f {1}{c^n}\big(b_n(x+d_n)u_n(x)u_n(y)-ca_nu_{n-1}(x)u_n(y)\big)
\\&\q-\f{1}{c^n}
\big(b_n(y+d_n)u_n(x)u_n(y)
-ca_nu_n(x)u_{n-1}(y)\big)
\\&=h_n(x,y)+\f{b_n}{c^n}(x-y)u_n(x)u_n(y).
\end{align*}
Thus,
$$(x-y)\sum_{n=0}^{p-1}\f{b_n}{c^n}u_n(x)u_n(y)
=\sum_{n=0}^{p-1}(h_{n+1}(x,y)
-h_n(x,y))=h_p(x,y)-h_0(x,y)=h_p(x,y).$$
This proves (i). By (i),
\begin{align*}\sum_{n=0}^{p-1}\f{b_n}{c^n}u_n(x)^2
&=\f{a_p}{c^{p-1}}\lim_{y\to x}\Big(u_{p-1}(x)\cdot
\f{u_p(x)-u_p(y)}{x-y}-u_p(x)\cdot\f{u_{p-1}(x)-u_{p-1}(y)}{x-y}\Big)
\\&=\f{a_p}{c^{p-1}}\big(u_{p-1}(x)u_p'(x)-u_p(x)u_{p-1}'(x)\big).
\end{align*}
The proof is now complete.
\vskip0.2cm
\par{\bf Definition 2.1} For $a,c\in\Bbb R$ with $c\not=0$ define $\{A_n(a,c;x)\}$ by
$A_0(a,c;x)=1$, $A_1(a,c;x)=x$ and $$(n+1)^3A_{n+1}(a,c;x)=(2n+1)(an(n+1)+x)A_n(a,c;x)
-cn^3A_{n-1}(a,c;x)\ (n\ge 1).$$
\par{\bf Corollary 2.1} Suppose that $p\in\Bbb Z^+$ and $a,c\in\Bbb R$ with $c\not=0$.
\par $(\t{\rm i})$ For $x\not=y$ we have
\begin{align*}&\sum_{n=0}^{p-1}\f{2n+1}{c^n}A_n(a,c;x)A_n(a,c;y)
\\&=\f{p^3}{c^{p-1}(x-y)}\big(A_p(a,c;x)A_{p-1}(a,c;y)
-A_{p-1}(a,c;x)A_p(a,c;y)\big).\end{align*}
\par $(\t{\rm ii})$
We have
\begin{align*}
&\sum_{n=0}^{p-1}\f{2n+1}{c^n}A_n(a,c;x)^2
\\&=\f{p^3}{c^{p-1}}
\Big(A_{p-1}(a,c;x)\f d{dx}A_p(a,c;x)-A_p(a,c;x)\f d{dx}A_{p-1}(a,c;x)\Big).
\end{align*}

\par{\it Proof.} Putting $a_n=n^3,\ b_n=2n+1$, $d_n=an(n+1)$ and $u_n(x)=A_n(a,c;x)$ in Theorem 2.1 gives the result.
\vskip0.2cm
\par{\bf Definition 2.2} For $a,c\in\Bbb R$ with $c\not=0$ define $\{A'_n(a,c;x)\}$ by
$A'_0(a,c;x)=1$, $A'_1(a,c;x)=x$ and $$(n+1)^2A'_{n+1}(a,c;x)=(an(n+1)+x)A'_n(a,c;x)
-cn^2A'_{n-1}(a,c;x)\ (n\ge 1).$$

\par{\bf Corollary 2.2} Suppose that $p\in\Bbb Z^+$ and $a,c\in\Bbb R$ with $c\not=0$.
\par $(\t{\rm i})$ For $x\not=y$ we have
\begin{align*}&\sum_{n=0}^{p-1}\f{1}{c^n}A'_n(a,c;x)A'_n(a,c;y)
\\&=\f{p^2}{c^{p-1}(x-y)}\big(A'_p(a,c;x)A'_{p-1}(a,c;y)
-A'_{p-1}(a,c;x)A'_p(a,c;y)\big).\end{align*}
\par $(\t{\rm ii})$
We have
\begin{align*}
\sum_{n=0}^{p-1}\f{1}{c^n}A'_n(a,c;x)^2
=\f{p^2}{c^{p-1}}
\Big(A'_{p-1}(a,c;x)\f d{dx}A'_p(a,c;x)-A'_p(a,c;x)\f d{dx}A'_{p-1}(a,c;x)\Big).
\end{align*}

\par{\it Proof.} Putting $a_n=n^2,\ b_n=1$, $d_n=an(n+1)$ and $u_n(x)=A'_n(a,c;x)$ in Theorem 2.1 gives the result.
\vskip0.2cm
\par{\bf Theorem 2.2} {\sl Suppose that $a_0=0$, $c\not=0$ and
$$a_{n+1}u_{n+1}(x)=b_n(x+d_n)u_n(x)-ca_nu_{n-1}(x)\ (n=1,2,3,\ldots).$$
For $p=1,2,3,\ldots$ we have}
$$\sum_{n=0}^{p-1}\f{b_n(x+y+2d_n)}{(-c)^n}u_n(x)u_n(y)
=\f{a_p}{(-c)^{p-1}}\big(u_p(x)u_{p-1}(y)+u_{p-1}(x)u_p(y)\big)$$
and
$$\sum_{n=0}^{p-1}\big(x-y+(-1)^{p-1-n}(x+y+2d_n)\big)\f{b_n}{c^n}
u_n(x)u_n(y)=\f{2a_p}{c^{p-1}}u_p(x)u_{p-1}(y).$$
\par{\it Proof.}
Set $$h_n(x,y)=\f{a_n}{(-c)^{n-1}}
\big(u_n(x)u_{n-1}(y)+u_{n-1}(x)u_n(y)\big).$$
Then
\begin{align*}
h_{n+1}(x,y)&=\f 1{(-c)^n}\big(a_{n+1}u_{n+1}(x)u_n(y)+ a_{n+1}u_{n+1}(y)u_n(x)\big)
\\&=\f 1{(-c)^n}\big(b_n(x+d_n)u_n(x)u_n(y)-ca_nu_{n-1}(x)u_n(y)
\\&\q+b_n(y+d_n)u_n(x)u_n(y)-ca_nu_n(x)u_{n-1}(y)\big)
\\&=\f{b_n(x+y+2d_n)}{(-c)^n}u_n(x)u_n(y)+h_n(x,y).
\end{align*}
Thus,
$$\sum_{n=0}^{p-1}\f{b_n(x+y+2d_n)}{(-c)^n}u_n(x)u_n(y)
=\sum_{n=0}^{p-1}(h_{n+1}(x,y)-h_n(x,y))
=h_p(x,y)-h_0(x,y)=h_p(x,y).$$
This proves the first identity. Combining this identity with Theorem 2.1 gives
\begin{align*}&\sum_{n=0}^{p-1}\big(x-y+(-1)^{p-1-n}(x+y+2d_n)\big)\f{b_n}{c^n}
\\&=\f{a_p}{c^{p-1}}(u_p(x)u_{p-1}(y)-u_{p-1}(x)u_p(y))
+\f {a_p}{c^{p-1}}(u_p(x)u_{p-1}(y)+u_{p-1}(x)u_p(y))
\\&=\f{2a_p}{c^{p-1}}u_p(x)u_{p-1}(y),\end{align*}
which concludes the proof.
\vskip0.2cm
\par{\bf Corollary 2.3} Let $a,c\in\Bbb R$ with $ac\not=0$. For $p\in\Bbb Z^+$ we have
\begin{align*}
&\sum_{n=0}^{p-1}
\f{(2n+1)^3}{(-c)^n}A_n(a,c;x)A_n(a,c;y)
\\&=\Big(1-\f 2a(x+y)\Big)\sum_{n=0}^{p-1}\f{2n+1}{(-c)^n}
A_n(a,c;x)A_n(a,c;y)
\\&\q+\f{2p^3}{a(-c)^{p-1}}\big(A_p(a,c;x)A_{p-1}(a,c;y)
+A_{p-1}(a,c;x)A_p(a,c;y)\big).
\end{align*}
\par{\it Proof.} Taking $a_n=n^3$, $b_n=2n+1$, $d_n=an(n+1)$ and $u_n=A_n(a,c;x)$ in Theorem 2.2 gives
\begin{align*}&\sum_{n=0}^{p-1}\f{(2n+1)(x+y+2an(n+1))}{(-c)^n}A_n(a,c;x)A_n(a,c;y)
\\&=\f{p^3}{(-c)^{p-1}}\big(A_p(a,c;x)A_{p-1}(a,c;y)
+A_{p-1}(a,c;x)A_p(a,c;y)\big).\end{align*}
Note that
$$(2n+1)(x+y+2an(n+1))=\f a2(2n+1)^3+\Big(x+y-\f a2\Big)(2n+1).$$
We then deduce the result.
\vskip0.2cm
\par{\bf Corollary 2.4} Let $a,c\in\Bbb R$ with $ac\not=0$. For $p\in\Bbb Z^+$ we have
\begin{align*}
&\sum_{n=0}^{p-1}
\f{(2n+1)^2}{(-c)^n}A'_n(a,c;x)A'_n(a,c;y)
\\&=\Big(1-\f 2a(x+y)\Big)\sum_{n=0}^{p-1}\f{1}{(-c)^n}
A'_n(a,c;x)A'_n(a,c;y)
\\&\q+\f{2p^2}{a(-c)^{p-1}}\big(A'_p(a,c;x)A'_{p-1}(a,c;y)
+A'_{p-1}(a,c;x)A'_p(a,c;y)\big).
\end{align*}
\par{\it Proof.} Taking $a_n=n^2$, $b_n=1$, $d_n=an(n+1)$ and $u_n=A'_n(a,c;x)$ in Theorem 2.2 and noting that $2n(n+1)=\f 12((2n+1)^2-1)$ yields the result.
\vskip0.2cm
\par{\bf Corollary 2.5} For $p=1,2,3,\ldots$ we have
$$\sum_{n=0}^{p-1}(x-y+(-1)^{p-1-n}(x+y))(2n+1)P_n(x)P_n(y)
=2pP_p(x)P_{p-1}(y),$$
where  $\{P_n(x)\}$ are the Legendre polynomials given by
$$P_0(x)=1,\  P_1(x)=x\ \t{and}\  (n+1)P_{n+1}(x)=(2n+1)xP_n(x)-nP_{n-1}(x)\ (n\ge 1).$$
\par{\it Proof.} Putting $a_n=n$, $b_n=2n+1$, $d_n=0$, $c=1$ and $u_n=P_n(x)$ in Theorem 2.2 gives the result.

\section*{3. Identities for $v_n(x)$ and $g_n(x)$}
\par\q For given positive real numbers $a,b,c,d$, the Wilson polynomial $W_n(x;a,b,c,d)$ is defined by
$$W_n(x;a,b,c,d)=(a+b)_n(a+c)_n(a+d)_n
\sum_{k=0}^n\f{(-n)_k(a+b+c+d+n-1)_k(a-i\sqrt x)_k(a+i\sqrt x)_k}{(a+b)_k(a+c)_k(a+d)_k\cdot k!},$$
where $(t)_0=1$ and $(t)_k=t(t+1)\cdots (t+k-1)$ for $k\ge 1$.
Since
$$(a+i\sqrt x)_k(a-i\sqrt x)_k=\prod_{j=0}^{k-1}(a+i\sqrt x+j)(a-i\sqrt x+j)=\prod_{j=0}^{k-1}\big(x+(a+j)^2\big)$$
and $(t)_k=\b{t-1+k}kk!=(-1)^k\b{-t}k k!$ we see that
\begin{align*}W_n(x;a,b,c,d)=&n!^3
\cdot (-1)^n\b{-a-b}n\b{-a-c}n\b{-a-d}n
\\&\times\sum_{k=0}^n\f{\b nk\b{n+a+b+c+d-2+k}k}
{\b{-a-b}k\b{-a-c}k\b{-a-d}k}\cdot
\frac{\prod_{j=0}^{k-1}(x+(a+j)^2)}{k!^2}.\tag{3.1}\end{align*}
It is well known (see [3]) that $\{W_n(x;a,b,c,d)\}$ are orthogonal polynomials on the half-line $(0,\infty)$ with respect to the weight function
$$w(x)=\f 1{2\sqrt x}\Big|\f{\Gamma(a+i\sqrt x)\Gamma(b+i\sqrt x)\Gamma(c+i\sqrt x)\Gamma(d+i\sqrt x)}{\Gamma(2i\sqrt x)}\Big|^2,$$
where $\Gamma(t)$ is the Gamma function. It is also known (see for example [3]) that $\{W_n(x;a,b,c,d)\}$ satisfy the recurrence relation:
\begin{align*} &W_{n+1}(x;a,b,c,d)
\\&=(A_n+C_n-a^2-x)\f{(2n-1+a+b+c+d)(2n+a+b+c+d)}
{n-1+a+b+c+d}W_n(x;a,b,c,d)
\\&\q-A_{n-1}C_n\f{\prod_{j=0}^3(2n-j+a+b+c+d)}
{(n-1+a+b+c+d)(n-2+a+b+c+d)}
 W_{n-1}(x;a,b,c,d)\ (n\ge 1),\tag{3.2}
\end{align*}
where
$$A_n=\f{(n-1+a+b+c+d)(n+a+b)(n+a+c)(n+a+d)}{(2n-1+a+b+c+d)(2n+a+b+c+d)}$$
and
$$C_n=\f{n(n-1+b+c)(n-1+b+d)(n-1+c+d)}{(2n-2+a+b+c+d)(2n-1+a+b+c+d)}.$$
\par{\bf Theorem 3.1} For $n=0,1,2,\ldots$ we have
$$v_n(x)=\f{W_n(\f{1-2x}4;\f 12,\f 12,\f 12,\f 12)}{n!^3}=\sum_{k=0}^n\f{\prod_{r=1}^k((2n+1)^2-(2r-1)^2)(2x-1-(2r-1)^2)}
{16^k\cdot k!^4}$$
and so
$$v_n(2m(m+1)+1)
=\sum_{k=0}^{\t{min}\{m,n\}}\f{\prod_{r=1}^k((2n+1)^2-(2r-1)^2)((2m+1)^2-(2r-1)^2)}
{16^k\cdot k!^4}$$
for $m=0,1,2,\ldots$
\vskip0.2cm
\par{\it Proof.} Note that
\begin{align*}
\b xk\b{x+k}k&=\f{x(x+1)\cdot (x-1)(x+2)\cdots (x-(k-1))(x+k)}{k!^2}
\\&=\f{x(x+1)\cdot (x(x+1)-2)\cdots(x(x+1)-k(k-1))}{k!^2}
\\&=\f{4x(x+1)\cdot(4x(x+1)-8)\cdots(4x(x+1)-4k(k-1))}
{4^k\cdot k!^2}.\end{align*}
We obtain
$$\b xk\b{x+k}k=\f{((2x+1)^2-1^2)((2x+1)^2-3^2)\cdots((2x+1)^2-(2k-1)^2)}
{4^k\cdot k!^2}.\eqno{(3.3)}$$
Hence,
\begin{align*}V_n(x)&=\sum_{k=0}^n\b nk\b{n+k}k\b xk\b{x+k}k
\\&=\sum_{k=0}^n
\f{\prod_{r=1}^k((2n+1)^2-(2r-1)^2)((2x+1)^2-(2r-1)^2)}{16^k\cdot k!^4}.\end{align*}
Since $V_n(x)=v_n(2x^2+2x+1)$, we see that
\begin{align*}
v_n(x)&=V_n\Big(\f{-1\pm \sqrt{2x-1}}2\Big)
=\sum_{k=0}^n
\f{\prod_{r=1}^k((2n+1)^2-(2r-1)^2)(2x-1-(2r-1)^2)}{16^k\cdot k!^4}.\end{align*}
Thus, for $m=0,1,2,\ldots$ we have
\begin{align*}
v_n(2m(m+1)+1)&=\sum_{k=0}^n
\f{\prod_{r=1}^k((2n+1)^2-(2r-1)^2)((2m+1)^2-(2r-1)^2)}{16^k\cdot k!^4}
\\&=\sum_{k=0}^{\t{min\{m,n\}}}
\f{\prod_{r=1}^k((2n+1)^2-(2r-1)^2)((2m+1)^2-(2r-1)^2)}{16^k\cdot k!^4}.\end{align*}
On the other hand, putting $a=b=c=d=\f 12$ in (3.1) and applying the above and the fact $\b{-1}k=(-1)^k$ gives
\begin{align*}
W_n\Big(\f{1-2x}4;\f 12,\f 12,\f 12,\f 12\Big)
&=n!^3\sum_{k=0}^n\b nk\b{n+k}k(-1)^k\cdot\f{\prod_{j=0}^{k-1}
(\f{1-2x}4+(\f{2j+1}2)^2)}{k!^2}
\\&=n!^3\sum_{k=0}^n\f{\prod_{r=1}^k((2n+1)^2-(2r-1)^2)(2x-1-(2r-1)^2)}
{16^k\cdot k!^4}
\\&=v_n(x),\end{align*}
which completes the proof.
\vskip0.2cm
By Corollary 2.1,
\begin{align*}&\sum_{n=0}^{p-1}(2n+1)v_n(x)v_n(y)
=\f{p^3}{x-y}(v_p(x)v_{p-1}(y)-v_{p-1}(x)v_p(y)),\tag{3.4}
\\&\sum_{n=0}^{p-1}(2n+1)v_n(x)^2=p^3(v_{p-1}(x)v'_p(x)-v_p(x)v'_{p-1}(x))
.\tag{3.5}\end{align*}
From Corollary 2.3,
\begin{align*}\sum_{n=0}^{p-1}(2n+1)^3(-1)^nv_n(x)v_n(y)
&=(1-2(x+y))\sum_{n=0}^{p-1}(2n+1)(-1)^nv_n(x)v_n(y)
\\&\q+2(-1)^{p-1}p^3(v_p(x)v_{p-1}(y)+v_{p-1}(x)v_p(y)).\tag{3.6}
\end{align*}
Since $V_n(x)=v_n(2x^2+2x+1)$, from (3.4) we have
\begin{align*}
&\sum_{n=0}^{p-1}(2n+1)V_n(x)V_n(y)
\\&=\f{p^3}{2(x-y)(x+y+1)}\big(V_p(x)V_{p-1}(y)-V_{p-1}(x)V_p(y)\big)
\\&=\f{p^3}{2(x+y+1)}\Big(V_{p-1}(x)\f{V_p(x)-V_p(y)}{x-y}
-V_p(x)\f{V_{p-1}(x)-V_{p-1}(y)}{x-y}\Big).
\end{align*}
Let $y\to x$. We get
$$\sum_{n=0}^{p-1}(2n+1)V_n(x)^2=\f{p^3}{2(2x+1)}
\big(V_{p-1}(x)V'_p(x)-V_p(x)V'_{p-1}(x)\big).\eqno{(3.7)}$$
\par{\bf Theorem 3.2} Let $p\in\Bbb Z^+$.
 \par $(\t{\rm i})$ For $x\not=1$ we have
$$\sum_{n=0}^{p-1}(2n+1)v_n(x)
=\f{p^3}{x-1}(v_{p}(x)-v_{p-1}(x)).$$
\par $(\t{\rm ii})$ For $x\not=1,5$ we have
\begin{align*}&
\sum_{n=0}^{p-1}(2n+1)^3v_n(x)
\\&=\f{p^3}{(x-1)(x-5)}\big(((2p-1)^2(x-1)+4)v_p(x)
-((2p+1)^2(x-1)+4)v_{p-1}(x)\big)
.\end{align*}
\par $(\t{\rm iii})$ For $x\not=1,5,13$ we have
\begin{align*}&
\sum_{n=0}^{p-1}(2n+1)^5v_n(x)
\\&=\f{p^3}{(x-1)(x-5)(x-13)}\big(((2p-1)^4(x-1)(x-5)
+48(2p-1)^2(x-1)+20x+92)v_p(x)
\\&\q-((2p+1)^4(x-1)(x-5)
+48(2p+1)^2(x-1)+20x+92)v_{p-1}(x)\big).
\end{align*}
\par{\it Proof.} Since $(n+1)^3=(2n+1)(n(n+1)+1)-n^3$, we see that $v_n(1)=1$. Thus, taking $y=1$ in (3.4) gives (i).
 Observe that
$$(n+1)^3(2(n+1)
(n+2)+1)=(2n+1)(n(n+1)+5)(2n(n+1)+1)-n^3(2n(n-1)+1).$$
We have
$$v_n(5)=2n(n+1)+1=\f{(2n+1)^2+1}2.\eqno{(3.8)}$$
Taking $y=5$ in (3.4) and then applying (3.8) yields
\begin{align*}\sum_{n=0}^{p-1}\f{(2n+1)^3+2n+1}2v_n(x)
=\f{p^3}{x-5}\big(v_p(x)(2(p-1)p+1)
-v_{p-1}(x)(2p(p+1)+1)\big).\end{align*}
This together with (i) gives
\begin{align*}
\sum_{n=0}^{p-1}(2n+1)^3v_n(x)
&=2\sum_{n=0}^{p-1}\f{(2n+1)^3+2n+1}2v_n(x)-\sum_{n=0}^{p-1}(2n+1)
v_n(x)
\\&=\f{2p^3}{x-5}\big(v_p(x)(2(p-1)p+1)
-v_{p-1}(x)(2p(p+1)+1)
\big)
\\&\q-\f{p^3}{x-1}(v_{p}(x)-v_{p-1}(x)),\end{align*}
which yields (ii).
\par Since $13=2\cdot 2\cdot 3+1$, from Theorem 3.1 we see that
$$v_n(13)=\f{3(2n+1)^4+18(2n+1)^2+11}{32}.\eqno{(3.9)}$$
Taking $y=13$ in (3.4) and then applying (3.9) yields
\begin{align*}&\sum_{n=0}^{p-1}(2n+1)v_n(x)
\f{3(2n+1)^4+18(2n+1)^2+11}{32}
\\&=\f{p^3}{x-13}\Big(v_p(x)\f{3(2p-1)^4+18(2p-1)^2+11}{32}
-v_{p-1}(x)\f{3(2p+1)^4+18(2p+1)^2+11}{32}\Big).\end{align*}
That is,
\begin{align*}3\sum_{n=0}^{p-1}(2n+1)^5v_n(x)
&=-18\sum_{n=0}^{p-1}(2n+1)^3v_n(x)-11\sum_{n=0}^{p-1}(2n+1)
v_n(x)
\\&\q+\f{p^3}{x-13}\big((3(2p-1)^4+18(2p-1)^2+11)v_p(x)
\\&\q-(3(2p+1)^4
+18(2p+1)^2+11)v_{p-1}(x)\big).
\end{align*}
This together with (i) and (ii) yields (iii).
\vskip0.2cm
\par{\bf Remark 3.1} In a similar way, one may deduce the identities for $\sum_{n=0}^{p-1}(2n+1)^kv_n(x)$ for $k=7,9,11,\ldots$
\vskip0.2cm
\par{\bf Theorem 3.3} For $n=0,1,2,\ldots$ we have
$$g_n(x)=\sum_{k=0}^n\b nk\f{\prod_{r=1}^k(2x-1-(2r-1)^2)}{4^k\cdot k!^2}$$ and so
$$g_n(2m^2+2m+1)=\sum_{k=0}^{\t{min}\{m,n\}}\b nk\f{\prod_{r=1}^k((2m+1)^2-(2r-1)^2)}{k!^2}.$$
\par{\it Proof.} By (3.3),
$$G_n(x)=\sum_{k=0}^n\b nk\b xk\b{x+k}k
=\sum_{k=0}^n\b nk\f{\prod_{r=1}^k((2x+1)^2-(2r-1)^2)}{4^k\cdot k!^2}.$$
Since $G_n(x)=g_n(2x^2+2x+1)$ we see that
$$g_n(x)=G_n\Big(\f{-1\pm \sqrt{2x-1}}2\Big)
=\sum_{k=0}^n\b nk\f{\prod_{r=1}^k(2x-1-(2r-1)^2)}{4^k\cdot k!^2}$$
It then follows that for $m\in\Bbb Z^+$,
$$g_n(2m^2+2m+1)=\sum_{k=0}^{\t{min}\{m,n\}}\b nk\f{\prod_{r=1}^k((2m+1)^2-(2r-1)^2)}{4^k\cdot k!^2},$$
which completes the proof.
\vskip0.2cm
By Corollary 2.2,
\begin{align*}
&\sum_{n=0}^{p-1}g_n(x)g_n(y)
=\f{2p^2}{x-y}(g_p(x)g_{p-1}(y)-g_{p-1}(x)g_p(y)),\tag{3.10}
\\&\sum_{n=0}^{p-1}g_n(x)^2
=2p^2(g_{p-1}(x)g'_p(x)-g_p(x)g'_{p-1}(x)).
\tag{3.11}
\end{align*}
From Corollary 2.4,
\begin{align*}&\sum_{n=0}^{p-1}(2n+1)^2(-1)^ng_n(x)g_n(y)
\\&=-\f{x+y}2\sum_{n=0}^{p-1}(-1)^ng_n(x)g_n(y)
+(-1)^{p-1}p^2(g_p(x)g_{p-1}(y)+g_{p-1}(x)g_p(y)).
\tag {3.12}\end{align*}
Since $G_n(x)=g_n(2x^2+2x+1)$, from (3.11) or (3.12) we deduce that
$$\sum_{n=0}^{p-1}G_n(x)^2=\f{p^2}{2x+1}\big(G'_p(x)G_{p-1}(x)
-G_p(x)G'_{p-1}(x)\big).\eqno{(3.13)}$$

\vskip0.2cm

\par{\bf Theorem 3.4} Let $p\in\Bbb Z^+$.
\par $(\t{\rm i})$ For $x\not=1$,
$$\sum_{n=0}^{p-1}g_n(x)
=\f{2p^2}{x-1}(g_p(x)-g_{p-1}(x)).$$
\par $(\t{\rm ii})$ For $x\not=1,5$,
$$\sum_{n=0}^{p-1}ng_n(x)
=\f{2p^2}{(x-1)(x-5)}\big((p(x-1)-(x-3))g_p(x)-(p(x-1)+2)
g_{p-1}(x)\big).$$
\par $(\t{\rm iii})$ For $x\not=1,5,13$,
\begin{align*}&\sum_{n=0}^{p-1}n^2g_n(x)
\\&=\f{2p^2}{(x-1)(x-5)(x-13)}
\big(\big((x-1)(x-5)p^2-2(x-1)(x-9)p+x^2-12x+19
\big)g_p(x)
\\&\q-\big((x-1)(x-5)p^2+8(x-1)p+2x+6\big)g_{p-1}(x)\big).
\end{align*}
\par{\it Proof.} Since $g_n(1)=1$, taking $y=1$ in (3.10) yields (i).
From Theorem 3.3,
$g_n(5)=2n+1,\ g_n(13)=3n(n+1)+1.$
Now, taking $y=5,13$ in (3.10) yields
\begin{align*}&\sum_{n=0}^{p-1}(2n+1)g_n(x)
=\f{2p^2}{x-5}\big((2p-1)g_p(x)-(2p+1)g_{p-1}(x)\big)\ (x\not=5),\tag{3.14}
\\&\sum_{n=0}^{p-1}(3n(n+1)+1)g_n(x)
\\&\q=\f{2p^2}{x-13}\big((3p(p-1)+1)g_p(x)-(3p(p+1)+1)g_{p-1}(x)\big)\ (x\not=13).\tag{3.15}\end{align*}
Combining (3.14) with (i) yields (ii). By (3.14), (3.15) and (i), for $x\not=1,5,13$,
\begin{align*}&3\sum_{n=0}^{p-1}n^2g_n(x)
\\&=\f 12\sum_{n=0}^{p-1}g_n(x)-\f 32\sum_{n=0}^{p-1}(2n+1)g_n(x)
+\f{2p^2}{x-13}\big((3p(p-1)+1)g_p(x)-(3p(p+1)+1)g_{p-1}(x)\big)
\\&=\f{p^2}{x-1}(g_p(x)-g_{p-1}(x))-\f{3p^2}{x-5}
\big((2p-1)g_p(x)-(2p+1)g_{p-1}(x)\big)
\\&\q+\f{2p^2}{x-13}\big((3p(p-1)+1)g_p(x)-(3p(p+1)+1)g_{p-1}(x)\big),
\end{align*}
which yields (iii).
\vskip0.2cm
\par{\bf Remark 3.2} In a similar way, one may deduce the identities for $\sum_{n=0}^{p-1}n^kg_n(x)$ for $k=3,4,5,\ldots$
\vskip0.2cm
\par{\bf Theorem 3.5} {\sl We have
$$v_n(x)=\sum_{k=0}^n\b nk\b{n+k}k(-1)^{n-k}g_k(x)\ (n=0,1,2,\ldots).$$
Moreover, for any odd prime $p$ and $m,x\in\Bbb Z_p$ with $m\not\e 0\mod p$ and $\ls{2x-1}p=1$,}
$$\sum_{k=0}^{p-1}\f{v_k(x)}{m^k}\e \Big(
\sum_{k=0}^{p-1}\f{g_k(x)}{m^k}\Big)^2
\mod p.$$
\par{\it Proof.}
Since $V_n(x)=v_n(2x^2+2x+1)$ and $G_n(x)=g_n(2x^2+2x+1)$, we see that
$$v_n(x)=V_n\Big(\f{-1\pm \sqrt{2x-1}}2\Big)\q\t{and}\q
g_n(x)=G_n\Big(\f{-1\pm \sqrt{2x-1}}2\Big).$$
Now the results follow from (1.1) and (1.2).
\vskip0.2cm
\par{\bf Theorem 3.6} {\sl Let $p$ be an odd prime and $x\in\Bbb Z_p$. Suppose that $\ls{2x-1}p=1$ and so $2x-1\e (2x_0-1)^2\mod p$ for some $x_0\in\{1,2,\ldots,\f{p-1}2\}$. For $n=0,1,\ldots,p-1$ we have}
$$g_n(x)\e G_n(x_0-1)\mod p\q\t{\sl and}\q
v_n(x)\e V_n(x_0-1)\mod p.$$
\par{\it Proof.} From Theorem 3.3,
 \begin{align*}
 g_n(x)&=\sum_{k=0}^n\b nk\f{\prod_{r=1}^k(2x-1-(2r-1)^2)}{4^k\cdot k!^2}
\\&\e 1+\sum_{k=1}^n\b nk\f{\prod_{r=1}^k((2x_0-1)^2-(2r-1)^2)}{4^k\cdot k!^2}
\\&=1+\sum_{k=1}^n\b nk\f{\prod_{r=1}^k(x_0-r)(x_0-1+r)}{k!^2}
=\sum_{k=0}^n\b nk \b{x_0-1}k\b{x_0-1+k}k
\\&=G_n(x_0-1)\mod p.
\end{align*}
This proves (i). On the other hand,
\begin{align*}
 v_n(x)&=\sum_{k=0}^n
 \f{\prod_{r=1}^k((2n+1)^2-(2r-1)^2)(2x-1-(2r-1)^2)}{16^k\cdot k!^4}
\\&\e 1+\sum_{k=1}^n
\f{\prod_{r=1}^k((2n+1)^2-(2r-1)^2)((2x_0-1)^2-(2r-1)^2)}
{16^k\cdot k!^4}
\\&=1+\sum_{k=1}^n\f{\prod_{r=1}^k(n+1-r)(n+r)
(x_0-r)(x_0-1+r)}{k!^4}
\\&=\sum_{k=0}^n\b nk \b{n+k}k\b{x_0-1}k\b{x_0-1+k}k
\\&=V_n(x_0-1)\mod p.
\end{align*}
This proves (ii). The proof is now complete.

\section*{4. Congruences involving $v_n(x)$ and $V_n(x)$}

\par{\bf Lemma 4.1} {\sl Let $p$ be an odd prime and $x\in\Bbb Z_p$ with $\ls{2x-1}p=1$. Then $v_{p-1}(x),v'_{p-1}(x),v_p(x),$ $pv'_p(x)\in\Bbb Z_p$.
Moreover, }
\begin{align*}&v_{p-1}(x)\e 1-2pH_{x_0-1}\mod {p^2},\
v'_{p-1}(x)\e 0\mod p, \\&v_p(x)\e
1+\f{x_1(x_1+8x_0-4)}{8(2x-1)}
+p\f{x_1^2(x_1+12x_0-6)}{16(2x-1)^2}
+\Big(2+\f{x_1(x_1+8x_0-4)}
{4(2x-1)}\Big)pH_{x_0-1}
\mod {p^2},
\\& pv'_p(x) \e
\f{x_1+4x_0-2}{2(2x-1)}
\Big(1+p\Big(\f {x_1+4x_0-2}{4(2x-1)}-\f 1{x_1+4x_0-2}+2H_{x_0-1}\Big)\Big)
\mod {p^2},\end{align*}
where $x_0$ and $x_1$ are given by $x_0\in\{1,2,\ldots,\f{p-1}2\}$, $2x-1\e (2x_0-1)^2\mod p$ and $2x-1=(2x_0-1)^2+px_1$.
\vskip0.2cm
\par{\it Proof}. From Theorem 3.1,
$$v_n(x)=\sum_{k=0}^n\f{\prod_{r=1}^k((2n+1)^2-(2r-1)^2)(2x-1-(2r-1)^2)}
{16^k\cdot k!^4}.$$
Thus $v_{p-1}(x),v'_{p-1}(x)\in\Bbb Z_p$.
  From [12, Lemma 2.3], for any $a\in\Bbb Z_p$,
  $$\sum_{k=1}^{p-1}\f{\b ak\b{-1-a}k}k\e -2H_{\langle a\rangle_p}\mod p,$$ where $\langle a\rangle_p$ is the least non-negative residue of $a$ modulo $p$. Actually, using the summation package Sigma in Mathematica one can check the identity
  $$\sum_{k=1}^n\b nk\b{n+k}k\f{(-1)^k}k=-2H_n.\eqno{(4.1)}$$
   Thus,
   $$\sum_{k=1}^{p-1}\f{\b ak\b{-1-a}k}k
   =\sum_{k=1}^{p-1}\b ak\b{a+k}k\f{(-1)^k}k\e \sum_{k=1}^{\ap}
   \b{\ap}k\b{\ap+k}k\f{(-1)^k}k=-2H_{\ap}\mod p.$$
   Therefore,
   \begin{align*}&\sum_{k=1}^{p-1}\f{\prod_{r=1}^k((2p\pm 1)^2-(2r-1)^2)(2x-1-(2r-1)^2)}
{16^k\cdot k!^4}
\\&\e \pm 4p\sum_{k=1}^{p-1}
\f{\prod_{r=2}^k(1-(2r-1)^2)\cdot\prod_{r=1}^k((2x_0-1)^2-(2r-1)^2)}
{16^k\cdot k!^4}
\\&\e \pm 4p\sum_{k=1}^{x_0-1}
\f{\prod_{r=2}^k(-4r(r-1))\cdot\prod_{r=1}^k4(x_0-r)(x_0+r-1)}
{16^k\cdot k!^4}
\\&\e \pm 4p\sum_{k=1}^{x_0-1}
\f{(-4)^{k-1}\cdot k!(k-1)!\cdot 4^k(x_0-1)\cdots(x_0-k)\cdot x_0(x_0+1)\cdots(x_0+k-1)}
{16^k\cdot k!^4}
\\&\e \mp p\sum_{k=1}^{x_0-1}\f{(-1)^k}k\b{x_0-1+k}k\b{x_0-1}k
 =\pm 2pH_{x_0-1}\mod {p^2}.
\end{align*}
Hence
$$v_{p-1}(x)=1+\sum_{k=1}^{p-1}\f{\prod_{r=1}^k((2p- 1)^2-(2r-1)^2)(2x-1-(2r-1)^2)}
{16^k\cdot k!^4}\e 1-2pH_{x_0-1}\mod {p^2}.$$
Note that $$v_{p-1}(x)=1+4p(p-1)\sum_{k=1}^{p-1}\f{\prod_{r=2}^k((2p- 1)^2-(2r-1)^2)\prod_{r=1}^k(2x-1-(2r-1)^2)}
{16^k\cdot k!^4}.$$
We see that $v'_{p-1}(x)\in\Bbb Z_p$ and
\begin{align*}v'_{p-1}(x)&=4p(p-1)\sum_{k=1}^{p-1}\f{\prod_{r=2}^k((2p- 1)^2-(2r-1)^2)\sum_{s=1}^k2\prod_{1\le r\le k,r\not=s}(2x-1-(2r-1)^2)}
{16^k\cdot k!^4}
\\&\q\e 0\mod p.
\end{align*}
\par On the other hand,
\begin{align*}&\prod_{r=1}^p((2p+1)^2-(2r-1)^2)
\\&\e 4p(p+1)\cdot 8p\cdot \prod_{r=2}
^{p-1}4(p-r(r-1))
\\&\e 32p^2(p+1)\cdot 4^{p-2}\prod_{r=2}^{p-1}(-r(r-1))\Big(
1-p\sum_{s=2}^{p-1}\f 1{s(s-1)}\Big)
\\&=32p^2(p+1)\cdot (-4)^{p-2}\cdot (p-1)!(p-2)!\Big(1-p\sum_{s=2}^{p-1}\Big(\f 1{s-1}-\f 1s\Big)\Big)
\\&=-8p^2(p+1)\cdot 4^{p-1}\cdot \f{(p-1)!^2}{p-1}\Big(1-p\Big (1-\f 1{p-1}\Big)\Big)
\\&\e 8p^2\cdot 4^{p-1}\cdot (p-1)!^2
\mod {p^4}
\end{align*}
and
\begin{align*}
&\sum_{1\le s\le p, s\not=x_0,p+1-x_0}\f 1{(2x_0-1)^2-(2s-1)^2}
\\&=\f 1{4(2x_0-1)}\sum_{1\le s\le p, s\not=x_0,p+1-x_0}\Big(\f 1{x_0-s}+\f 1{x_0+s-1}\Big)
\\&=\f 1{4(2x_0-1)}\Big(\sum_{1\le s\le p,s\not=x_0}\f 1{x_0-s}-\f 1{x_0-(p+1-x_0)}+\sum_{1\le s\le p,s\not=p+1-x_0}\f 1{x_0+s-1}-\f 1{2x_0-1}\Big)
\\&\e \f 1{4(2x_0-1)}\Big(-\f 1{2x_0-1}-\f 1{2x_0-1}\Big)
= -\f 1{2(2x_0-1)^2}\e -\f 1{2(2x-1)}\mod p.\end{align*}
Thus,
\begin{align*}
&\prod_{1\le r\le p,\ r\not=x_0,p+1-x_0} (2x-1-(2r-1)^2)
\\&= \prod_{1\le r\le p,\ r\not=x_0,p+1-x_0} (px_1+(2x_0-1)^2-(2r-1)^2)
\\&=\prod_{1\le r\le p,\ r\not=x_0,p+1-x_0} ((2x_0-1)^2-(2r-1)^2) \Big(1+px_1\sum_{1\le s\le p,s\not=x_0,p+1-x_0}\f 1{(2x_0-1)^2-(2s-1)^2}\Big)
\\&\e \prod_{1\le r\le p,\ r\not=x_0,p+1-x_0}
4(x_0-r)(x_0+r-1)\cdot \Big(1-\f{px_1}{2(2x-1)}\Big)
\\&=4^{p-2}\cdot \f{\prod_{1\le r\le p,\ r\not=x_0} (x_0-r)}
{x_0-(p+1-x_0)}\cdot \f {\prod_{1\le r\le p,\ r\not=p+1-x_0}  (x_0+r-1)}{x_0+x_0-1}\cdot\Big(1-\f{px_1}{2(2x-1)}\Big)
\\&=4^{p-2}\cdot\f{(x_0-1)!(-1)^{p-x_0}(p-x_0)!\cdot x_0(x_0+1)\cdots(p-1)\cdot (p+1)\cdots(p+x_0-1)}{(2x_0-1)(2x_0-1-p)}
\\&\q\times\Big(1-\f{px_1}{2(2x-1)}\Big)
\\&\e 4^{p-2}(-1)^{x_0-1}\f{(p-1)!(p-x_0)!(x_0-1)!(1+pH_{x_0-1})}
{(2x_0-1)^2-p(2x_0-1)}\Big(1-\f{px_1}{2(2x-1)}\Big)
\\&=4^{p-2}\f{(p-1)!^2(1+pH_{x_0-1})}
{((2x_0-1)^2-p(2x_0-1))(-1)^{x_0-1}\b{p-1}{x_0-1}}\Big(1-\f{px_1}{2(2x-1)}\Big)
\\&\e 4^{p-2}\cdot (p-1)!^2\f{1+pH_{x_0-1}} {(2x_0-1)^2(1-\f p{2x_0-1})(1-pH_{x_0-1})}
\Big(1-\f{px_1}{2(2x-1)}\Big)
\\&\e 4^{p-2}\cdot (p-1)!^2\f{(1+\f p{2x_0-1})(1+pH_{x_0-1})^2} {(2x_0-1)^2}\Big(1-\f{px_1}{2(2x-1)}\Big)
\\&\e 4^{p-2}\cdot (p-1)!^2\f 1{2x-1-px_1}
\Big(1+p\Big(\f 1{2x_0-1}+2H_{x_0-1}-\f {x_1}{2(2x-1)}\Big)\Big)
\\&\e \f{4^{p-1}\cdot (p-1)!^2}{4(2x-1)}
\Big(1+p\Big(\f 1{2x_0-1}+2H_{x_0-1}+\f {x_1}{2(2x-1)}\Big)\Big)
\mod {p^2}
\end{align*}
and so
\begin{align*}&\prod_{r=1}^p(2x-1-(2r-1)^2)
\\&=(2x-1-(2x_0-1)^2)(2x-1-(2(p+1-x_0)-1)^2)
\prod_{1\le r\le p,\ r\not=x_0,p+1-x_0} (2x-1-(2r-1)^2)
\\&=px_1(px_1+4p(2x_0-1)-4p^2)\prod_{1\le r\le p,\ r\not=x_0,p+1-x_0} (2x-1-(2r-1)^2)
\\&\e p^2x_1(x_1+4(2x_0-1)-4p)\cdot \f{4^{p-1}\cdot (p-1)!^2}{4(2x-1)}
\\&\q\times\Big(1+p\Big(\f 1{2x_0-1}+2H_{x_0-1}+\f {x_1}{2(2x-1)}\Big)\Big)\mod {p^4}.
\end{align*}
Therefore,
\begin{align*}&v_p(x)\\&=1+\f{\prod_{r=1}^p((2p+1)^2-(2r-1)^2)(2x-1-(2r-1)^2)}{16^p
\cdot p!^4}
\\&\q+\sum_{k=1}^{p-1}\f{\prod_{r=1}^k((2p+1)^2-(2r-1)^2)(2x-1-(2r-1)^2)}
{16^k\cdot k!^4}
\\&\e 1+\f{8p^2\cdot 4^{p-1}\cdot (p-1)!^2\cdot p^2}{16^p\cdot p^4\cdot(p-1)!^4}\cdot x_1(x_1+4(2x_0-1)-4p)
\\&\q\times\f{4^{p-1}\cdot (p-1)!^2}{4(2x-1)}\Big(1+p\Big(\f 1{2x_0-1}+2H_{x_0-1}+\f {x_1}{2(2x-1)}\Big)\Big)+2pH_{x_0-1}
\\&=1+2pH_{x_0-1}+
\f{1}{8(2x-1)}\cdot x_1(x_1+4(2x_0-1)-4p)
\\&\q\times\Big(1+p\Big(\f 1{2x_0-1}+2H_{x_0-1}+\f {x_1}{2(2x-1)}\Big)\Big)
\\&\e 1+2pH_{x_0-1}-\f{px_1}{2(2x-1)}
\\&\q+\f {x_1(x_1+4(2x_0-1))}{8(2x-1)}\Big(1+p\Big(\f 1{2x_0-1}+2H_{x_0-1}+\f {x_1}{2(2x-1)}\Big)\Big)
\\&\e 1+\f{x_1(x_1+8x_0-4)}{8(2x-1)}
+p\f{x_1^2(x_1+12x_0-6)}{16(2x-1)^2}
+\Big(2+\f{x_1(x_1+8x_0-4)}
{4(2x-1)}\Big)pH_{x_0-1}
\mod {p^2} .\end{align*}
\par Set $f_p(t)=\prod_{r=1}^p(2t-1-(2r-1)^2)$. Then
$f'_p(t)=\sum_{s=1}^p2\prod_{1\le r\le p,r\not=s}\big(2t-1-(2r-1)^2\big).$
Hence,
\begin{align*}\f 12f'_p(x)&
=\prod_{1\le r\le p, r\not=x_0,p+1-x_0}(2x-1-(2r-1)^2)
\times\Big(2x-1-(2(p+1-x_0)-1)^2
\\&\q+2x-1-(2x_0-1)^2
+(2x-1-(2x_0-1)^2)
(2x-1-(2(p+1-x_0)-1)^2)
\\&\q\times\sum_{1\le s\le p, s\not=x_0,p+1-x_0}\f 1{2x-1-(2s-1)^2}\Big)
\\&=\prod_{1\le r\le p, r\not=x_0,p+1-x_0}(2x-1-(2r-1)^2)
\times\Big(px_1+(2x_0-1)^2-(2(p+1-x_0)-1)^2
\\&\q+px_1+px_1(px_1+(2x_0-1)^2-(2(p+1-x_0)-1)^2)
\\&\q\times\sum_{1\le s\le p, s\not=x_0,p+1-x_0}\f 1{px_1+(2x_0-1)^2-(2s-1)^2}\Big)
\\&=p\prod_{1\le r\le p, r\not=x_0,p+1-x_0}(2x-1-(2r-1)^2)
\times\Big(2x_1+4(2x_0-1)-4p
\\&\q+px_1(x_1+4(2x_0-1)-4p)
\sum_{1\le s\le p, s\not=x_0,p+1-x_0}\f 1{(2x_0-1)^2-(2s-1)^2}\Big)
\\&\e p\cdot \f{4^{p-1}\cdot (p-1)!^2}{4(2x-1)}
\Big(1+p\Big(\f 1{2x_0-1}+2H_{x_0-1}+\f {x_1}{2(2x-1)}\Big)\Big)
\\&\q\times \Big(2x_1+4(2x_0-1)-4p
+px_1(x_1+4(2x_0-1))\f {-1}{2(2x-1)}\Big)
\mod {p^3}.
\end{align*}
Since
\begin{align*}pv_p(x)&=p+p\f{\prod_{r=1}^p((2p+1)^2-(2r-1)^2)
\cdot f_p(x)}{16^p
\cdot p!^4}
\\&\q+p\cdot 4p(p+1)\sum_{k=1}^{p-1}\f{\prod_{r=2}^k((2p+1)^2-(2r-1)^2)
\prod_{r=1}^k(2x-1-(2r-1)^2)}
{16^k\cdot k!^4},\end{align*}
we see that
\begin{align*}pv'_p(x)&\e p\f{\prod_{r=1}^p((2p+1)^2-(2r-1)^2)
\cdot f'_p(x)}{16^p
\cdot p!^4}
\\&\e \f{p\cdot 8p^2\cdot 4^{p-1}\cdot (p-1)!^2}{16^p\cdot p!^4}\cdot 2p\cdot \f{4^{p-1}\cdot (p-1)!^2}{4(2x-1)}
\Big(1+p\Big(\f 1{2x_0-1}+2H_{x_0-1}+\f {x_1}{2(2x-1)}\Big)\Big)
\\&\q\times \Big(2x_1+4(2x_0-1)-4p
-px_1(x_1+4(2x_0-1))\f {1}{2(2x-1)}\Big)
\\&\e \f{x_1+4x_0-2}{2(2x-1)}
\Big(1+p\Big(\f {x_1+4x_0-2}{4(2x-1)}-\f 1{x_1+4x_0-2}+2H_{x_0-1}\Big)\Big)
\mod {p^2}.
\end{align*}
This completes the proof.
\vskip0.2cm
\par{\bf Theorem 4.1} {\sl Let $p$ be an odd prime, $x\in\Bbb Z_p$ and $\ls{2x-1}p=1$. Let $x_0$ and $x_1$ be given by $x_0\in\{1,2,\ldots,\f{p-1}2\}$, $2x-1\e (2x_0-1)^2\mod p$ and $2x-1=(2x_0-1)^2+px_1$.
\par $(\t{\rm i})$ If $x\not\e 1\mod p$, then
\begin{align*}\sum_{n=0}^{p-1}(2n+1)v_n(x)
&\e p^3\f{x_1(x_1+8x_0-4)}{8(x-1)(2x-1)}
+\f{p^4}{x-1}
\Big(\f{x_1^2(x_1+12x_0-6)}{16(2x-1)^2}
\\&\q+\Big(4+\f{x_1(x_1+8x_0-4)}
{4(2x-1)}\Big)H_{x_0-1}\Big)
\mod {p^5}.\end{align*}
\par $(\t{\rm ii})$ If $x\not\e 1,5\mod p$, then
\begin{align*}
&\sum_{n=0}^{p-1}(2n+1)^3v_n(x)\\&\e p^3\f{(x+3)x_1(x_1+8x_0-4)}{8(x-1)(x-5)(2x-1)}
+\f{p^4}{x-5}\Big(\f{(x+3)x_1^2(x_1+12x_0-6)}{16(x-1)(2x-1)^2}
\\&\q-8-\f{x_1(x_1+8x_0-4)}{2(2x-1)}
+\f{x+3}{x-1}\Big(4+\f{x_1(x_1+8x_0-4)}{4(2x-1)}\Big)H_{x_0-1}
\Big)
\mod {p^5}.\end{align*}
\par $(\t{\rm iii})$ If $x\not\e 1,5,13\mod p$, then
\begin{align*}&\sum_{n=0}^{p-1}(2n+1)^5v_n(x)
\\&\e p^3\f{(x^2+62x+49)x_1(x_1+8x_0-4)}{8(x-1)(x-5)(x-13)(2x-1)}
\\&\q+\f{p^4(x^2+62x+49)}{(x-1)(x-5)(x-13)}
\Big(\f{x_1^2(x_1+12x_0-6)}{16(2x-1)^2}+\Big(4+\f{x_1(x_1+8x_0-4)}{4(2x-1)}
\Big)H_{x_0-1}\Big)
\\&\q-\f{p^4(x+19)}{(x-5)(x-13)}
\Big(16+\f{x_1(x_1+8x_0-4)}{2x-1}\Big)
\mod {p^5}.\end{align*}}
\par{\it Proof.} This is immediate from Theorem 3.2 and Lemma 4.1.
\vskip0.2cm
\par{\bf Remark 4.1.} Let $p>3$ be a prime and $x\in\Bbb Z_p$. In a previous preprint, the author established a congruence for
$\sum_{n=0}^{p-1}(2n+1)V_n(x)$ modulo $p^5$. Then Mao and Yang [8]  determined
$\sum_{n=0}^{\infty}(2n+1)^3V_n(x)$ and  $\sum_{n=0}^{\infty}(2n+1)^5V_n(x)$ modulo $p^4$ by finding suitable combinatorial identities.
\vskip0.2cm

\par{\bf Corollary 4.1} {\sl Let $p$ be a prime greater than $7$. Then}
\begin{align*}
&\sum_{n=0}^{p-1}(2n+1)^3V_n\Big(\f 12\Big)=\sum_{n=0}^{p-1}(2n+1)^3v_n\Ls 52\e
\f{11}{15}p^3+\f{p^4}{5}(-32+44q_p(2))\mod {p^5},
\\&\sum_{n=0}^{p-1}(2n+1)^3V_n\Big(-\f 13\Big)
=\sum_{n=0}^{p-1}(2n+1)^3v_n\Ls 59
\e -\f 45p^3+\f{p^4}5(7-42q_p(3))\mod {p^5},
\\&\sum_{n=0}^{p-1}(2n+1)^3V_n\Big(-\f 14\Big)=\sum_{n=0}^{p-1}(2n+1)^3v_n\Ls 58\e
-\f{29}{35}p^3+\f{p^4}{35}(52-754q_p(2))\mod {p^5},
\\&\sum_{n=0}^{p-1}(2n+1)^5V_n\Big(\f 12\Big)=\sum_{n=0}^{p-1}(2n+1)^5v_n\Ls 52\e
-\f{841}{315}p^3+\f{p^4}{105}(2332-3364q_p(2))\mod {p^5},
\\&\sum_{n=0}^{p-1}(2n+1)^5V_n\Big(-\f 13\Big)=\sum_{n=0}^{p-1}(2n+1)^5v_n\Ls 59\e
\f{53}{35}p^3+\f{p^4}{10}(-44+159q_p(3))\mod {p^5},
\\&\sum_{n=0}^{p-1}(2n+1)^5V_n\Big(-\f 14\Big)\\&=\sum_{n=0}^{p-1}(2n+1)^5v_n\Ls 58\e
\f{5641}{3465}p^3+\f{p^4}{3465}(-16328+146666q_p(2))\mod {p^5},
\\&\sum_{n=0}^{p-1}(2n+1)^3V_n\Big(-\f 16\Big)
\\&=\sum_{n=0}^{p-1}(2n+1)^3v_n\Ls {13}{18}
\e -\f{67}{77}p^3+\Big(\f{124}{77}-\f{2077}{385}(4q_p(2)+3q_p(3))\Big)p^4
\mod {p^5},
\\&\sum_{n=0}^{p-1}(2n+1)^5V_n\Big(-\f 16\Big)
\\&=\sum_{n=0}^{p-1}(2n+1)^5v_n\Ls {13}{18}
\e \f{30553}{17017}p^3+\Big(-\f{88040}{17017}+\f{947143}{85085}
(4q_p(2)+3q_p(3))\Big)p^4\mod {p^5}.
\end{align*}
\par{\it Proof.} It is well known (see for example [11]) that
\begin{align*}&H_{\f{p-1}2}\e -2q_p(2)\mod p,\q H_{[\f p3]}\e -\f 32q_p(3)\mod p,
\\&
\ H_{[\f p4]}\e -3q_p(2)\mod p,\q H_{[\f p6]}\e -2q_p(2)-\f 32q_p(3)\mod p.\end{align*}
Thus, the results follow from Theorem 4.1 and some easy computations. For example, taking $x=\f 52$, $x_0=\f{p-1}2$ and $x_1=4-p$ in Theorem 4.1 and then applying the facts $V_n(t)=v_n(2t^2+2t+1)$ and $H_{\f{p-1}2}\e -2q_p(2)\mod p$ yields the first and fourth congruences.
\vskip0.2cm

\par{\bf Lemma 4.2 ([6,7])} {\sl Let $p$ be a prime greater than $3$. Then}
\begin{align*}
&V_p\e 8+40p^3B_{p-3}\mod {p^4},
\\&V_p^{(3)}\e 15+132p^3B_{p-3}\mod
{p^4},
\\&V_p^{(4)}\e 40+704p^3B_{p-3}\mod {p^4},
\\& V_p^{(6)}\e 312+
16120p^3B_{p-3}\mod {p^4},
\\&V_{p-1}\e 256^{p-1}-\f 32p^3B_{p-3}\mod {p^4},
\\& V_{p-1}^{(3)}\e 729^{p-1}-\f{92}{27}p^3B_{p-3}\mod {p^4},
\\& V_{p-1}^{(4)}\e 4096^{p-1}-\f{17}2p^3B_{p-3}\mod {p^4}, \q
\\&V_{p-1}^{(6)}\e 186624^{p-1}-\f{1705}{54}p^3B_{p-3} \mod {p^4}.\end{align*}

\par{\bf Theorem 4.2} {\sl Let $p$ be a prime greater than $3$. Then}
\begin{align*}
&\sum_{n=0}^{p-1}(2n+1)^3\f{V_n}{16^n} \\&\e\f{p^3}{9}\Big((7+4p-4p^2)\cdot 16^{-(p-1)}-2(7-4p-4p^2)\cdot 16^{p-1}\Big)+\f{56}9p^6B_{p-3}\mod {p^7}
\end{align*}
and
\begin{align*}&\sum_{n=0}^{p-1}(2n+1)^5\f{V_n}{16^n}
\e \f{p^3}{75}\big((96p^3+56p^2-104p-107)\cdot 16^{-(p-1)}
\\&\qq\qq\q\qq\q+2(96p^3-56p^2-104p+107)\cdot 16^{p-1}\big)-\f{856}{75}p^6B_{p-3}\mod {p^7}.
\end{align*}

\par{\it Proof.} Since $V_n(t)=v_n(2t^2+2t+1)$, we see that
$\f{V_n}{16^n}=V_n(-\f 12)=v_n(\f 12)$.
Taking $x=\f 12$ in Theorem 3.2 and then applying Lemma 4.2 gives
\begin{align*}&\sum_{n=0}^{p-1}(2n+1)^3\f{V_n}{16^n}
=\sum_{n=0}^{p-1}(2n+1)^3v_n\Ls 12
\\&=\f{2p^3}9\Big((7+4p-4p^2)v_p\Ls 12-(7-4p-4p^2)v_{p-1}\Ls 12\Big)
\\&=\f{2p^3}9\Big((7+4p-4p^2)\f{V_p}{16^p}-(7-4p-4p^2)\f{V_{p-1}}{16^{p-1}}
\Big)
\\&\e \f{2p^3}9\Big((7+4p-4p^2)\f{8+40p^3B_{p-3}}{16^p}-(7-4p-4p^2)
\f{256^{p-1}-\f 32p^3B_{p-3}}{16^{p-1}}
\Big)
\\&\e\f{p^3}{9}\Big((7+4p-4p^2)\cdot 16^{-(p-1)}-2(7-4p-4p^2)\cdot 16^{p-1}\Big)+\f{56}9p^6B_{p-3}\mod {p^7}
\end{align*}
and
\begin{align*}\sum_{n=0}^{p-1}(2n+1)^5\f{V_n}{16^n}
&=\sum_{n=0}^{p-1}(2n+1)^5v_n\Ls 12
\\&=-\f{8p^3}{225}\Big(\big(\f 94(2p-1)^4-24(2p-1)^2+102\big)v_p\Ls 12
\\&\q-\big(\f 94(2p+1)^4-24(2p+1)^2+102\big)v_{p-1}\Ls 12\Big)
\\&=-\f{2p^3}{75}\Big(\big(3(2p-1)^4-32(2p-1)^2+136\big)\f{V_p}{16^p}
\\&\q-\big(3(2p+1)^4-32(2p+1)^2+136\big)\f{V_{p-1}}{16^{p-1}}
\\&\q=-\f{2p^3}{75}\Big(\big(-96p^3-56p^2+104p+107\big)\f{8+40p^3B_{p-3}}
{16^p}
\\&\q-\big(96p^3-56p^2-104p+107\big)\f {256^{p-1}-\f 32p^3B_{p-3}}{16^{p-1}}\Big)
\\&\e\f{p^3}{75}\big((96p^3+56p^2-104p-107)\cdot 16^{-(p-1)}
\\&\q+2(96p^3-56p^2-104p+107)\cdot 16^{p-1}\big)-\f{856}{75}p^6B_{p-3}\mod {p^7}.
\end{align*}
This proves the theorem.
\vskip0.2cm
\par{\bf Remark 4.2}
Since $V_n(x)=v_n(2x^2+2x+1)$ we see that
$$\ V_n\Big(-\f 13\Big)=v_n\Ls 59,\
V_n\Big(-\f 14\Big)=v_n\Ls 58,\ V_n\Big(-\f 16\Big)=v_n\Ls{13}{18}.$$
Similarly, using Theorem 3.2 and Lemma 4.2 one may deduce the congruences for $\sum_{n=0}^{p-1}(2n+1)^3V_n(-\f 1m)$ and
$\sum_{n=0}^{p-1}(2n+1)^5V_n(-\f 1m)$ modulo $p^7$ for $m=3,4,6$, where $p$ is a prime greater than $5$. The congruences for $\sum_{n=0}^{p-1}(2n+1)V_n(-\f 1m)$ $(m=2,3,4,6)$ modulo $p^7$ have been given by the author in [15, Theorem 22.47 (with $t=0$)].
\vskip0.2cm
\par{\bf Theorem 4.3} {\sl Let $p$ be an odd prime, $x\in\Bbb Z_p$, $\ls{2x-1}p=1$, $m\in\{0,1,2,\ldots\}$ and $p>2m$. Suppose that $x\not\e 2k(k+1)+1\mod p$ for $k=0,1,\ldots,m$. Then}
$$\sum_{n=0}^{p-1}(2n+1)^{2m+1}v_n(x)\e 0\mod {p^3}.$$
\par{\it Proof.} We prove the result by induction on $m$. For $m=0,1,2$ the result follows from Theorem 4.1. Now assume that $m\ge 3$ and
$\sum_{n=0}^{p-1}(2n+1)^{2r+1}v_n(x)\e 0\mod {p^3}$ for $r=0,1,\ldots,m-1$. By (3.4) and Lemma 4.1,
\begin{align*}&\sum_{n=0}^{p-1}(2n+1)v_n(2m(m+1)+1)v_n(x)
\\&=\f{p^3}{x-(2m(m+1)+1)}
(v_p(x)v_{p-1}(2m(m+1)+1)-v_{p-1}(x)v_p(2m(m+1)+1))
\\&\e 0
\mod {p^3}.\end{align*}
Since $p>2m$, from Theorem 3.1 we see that
\begin{align*}v_n(2m(m+1)+1)
&=\sum_{k=0}^{m}\f{\prod_{r=1}^k((2n+1)^2-(2r-1)^2)((2m+1)^2-(2r-1)^2)}
{16^k\cdot k!^4}
\\&=a_m(2n+1)^{2m}+a_{m-1}(2n+1)^{2(m-1)}
+a_1(2n+1)^2+a_0,\end{align*}
where $a_0,a_1,\ldots,a_m\in\Bbb Z_p$ and
$$a_m=\f{\prod_{r=1}^m((2m+1)^2-(2r-1)^2)}{16^m\cdot m!^4}
\not\e 0\mod p.$$
Hence,
$$a_m\sum_{n=0}^{p-1}(2n+1)^{2m+1}v_n(x)
=-\sum_{k=0}^{m-1}a_k\sum_{n=0}^{p-1}(2n+1)^{2k+1}v_n(x)
\e 0\mod {p^3}$$
by the inductive hypothesis. The proof is now complete.
\vskip0.2cm
\par{\bf Theorem 4.4} {\sl Let $p$ be an odd prime, $x\in\Bbb Z_p$ and $\ls{2x-1}p=1$. Then
$$\sum_{n=0}^{p-1}(2n+1)v_n(x)^2\e
\f{x_1+4x_0-2}{2(2x-1)}p^2+\f{x_1(x_1+8x_0-4)}{8(2x-1)^2}p^3
\mod {p^4},$$
where $x_0$ and $x_1$ are given by $x_0\in\{1,2,\ldots,\f{p-1}2\}$, $2x-1\e (2x_0-1)^2\mod p$ and $2x-1=(2x_0-1)^2+px_1$.}
\vskip0.2cm
\par{\it Proof.} From (3.3) and Lemma 4.1, \begin{align*}
\sum_{n=0}^{p-1}(2n+1)v_n(x)^2
&=p^3(v_{p-1}(x)v'_p(x)-v_p(x)v'_{p-1}(x))
\e p^2v_{p-1}(x)\cdot pv'_p(x)
\\&\e p^2(1-2pH_{x_0-1})\cdot\f{x_1+4x_0-2}{4x-2}
\\&\q\times\Big(1+p\Big(\f{x_1+4x_0-2}{4(2x-1)}-\f 1{x_1+4x_0-2}\Big)+2pH_{x_0-1}\Big)
\\&\e p^2\f{x_1+4x_0-2}{4x-2}\Big(1+p\Big(\f{x_1+4x_0-2}{4(2x-1)}-\f 1{x_1+4x_0-2}\Big)\Big)
\mod {p^4},\end{align*}
which yields the result.
\vskip0.2cm
\par{\bf Corollary 4.2} Let $p$ be an odd prime and $x\in\Bbb Z_p$ with $x\not\e -\f 12\mod p$. Then
$$\sum_{n=0}^{p-1}(2n+1)V_n(x)^2\e \f{1+2x'}{1+2x}p^2\mod
{p^4}, $$ where $x'=(x-\xp)/p$ and $\xp$ is the least non-negative residue of $x$ modulo $p$.
\vskip0.2cm
\par{\it Proof.} Observe that $V_n(x)=v_n(2x^2+2x+1)$ and
$$2(2x^2+2x+1)-1=(2x+1)^2\e (2(\xp+1)-1)^2
\e (2(p-\xp)-1)^2\mod p.$$
Set
$$ x_0=\begin{cases}\xp+1=x+1-px'&\t{if $\xp<\f{p-1}2$,}
\\p-\xp=(x'+1)p-x&\t{if $\xp\ge \f{p-1}2$}
\end{cases}$$
and
$$ x_1=\f{(2x+1)^2-(2x_0-1)^2}p=\begin{cases}
4x'(2x+1-px')&\t{if $\xp<\f{p-1}2$,}
\\ 4(x'+1)(2x+1-p(x'+1))&\t{if $\xp\ge \f{p-1}2$.}
\end{cases}$$
Then $x_0\in\{1,2,\ldots,\f{p-1}2\}$ and
$2(2x^2+2x+1)-1=(2x_0-1)^2+px_1$. From Theorem 4.4 we have
 \begin{align*}\sum_{n=0}^{p-1}(2n+1)V_n(x)^2&=
 \sum_{n=0}^{p-1}v_n(2x^2+2x+1)^2
 \\&\e
\f{x_1+4x_0-2}{2(2x+1)^2}p^2
+\f{x_1(x_1+8x_0-4)}{8(2x+1)^4}p^3
\mod {p^4}.\end{align*}
It is easy to see that
$$x_1+4x_0-2=2(2x+1)(2x'+1)-4px'(x'+1)$$
and
$$x_1+8x_0-4\e 16(2x+1)^2x'(x'+1)\mod p.$$
Thus,
\begin{align*}
\sum_{n=0}^{p-1}(2n+1)V_n(x)^2
&\e
\f{2(2x+1)(2x'+1)-4px'(x'+1)}{2(2x+1)^2}p^2
+\f{16(2x+1)^2x'(x'+1)}{8(2x+1)^4}p^3
\\&=\f{2x'+1}{2x+1}p^2\mod {p^4}.
\end{align*}
This proves the corollary.
\vskip0.2cm
\par{\bf Corollary 4.3} {\sl Let $p>3$ be a prime. Then
\begin{align*} &\sum_{n=0}^{p-1}(2n+1)V_n\Ls 12^2
=\sum_{n=0}^{p-1}(2n+1)v_n\Ls 52^2
\e 0\mod {p^4},
\\&\sum_{n=0}^{p-1}(2n+1)V_n\Big(-\f 13\Big)^2
=\sum_{n=0}^{p-1}(2n+1)v_n\Ls 59^2\e (-1)^{[\f p3]}p^2\mod {p^4},
\\&\sum_{n=0}^{p-1}(2n+1)V_n\Big(-\f 14\Big)^2
=\sum_{n=0}^{p-1}(2n+1)v_n\Ls 58^2\e (-1)^{\f{p-1}2}p^2\mod {p^4},
\\&\sum_{n=0}^{p-1}(2n+1)V_n\Big(-\f 16\Big)^2
=\sum_{n=0}^{p-1}(2n+1)v_n\Ls {13}{18}^2\e (-1)^{[\f p3]}p^2\mod {p^4}.
\end{align*}}
\par{\it Proof.} Since $V_n(x)=v_n(2x^2+2x+1)$, the result follows from Theorem 4.4. For example, putting $x=\f 58$, $x_0=\f 14(p+2+(-1)^{\f{p-1}2})$ and $x_1=-\f 14(p+2(-1)^{\f{p-1}2})$ in Theorem 4.4 yields the third congruence.

\vskip0.2cm

\par{\bf Remark 4.3}
In [16] the author conjectured that for any prime $p>3$,
\begin{align*} &\sum_{n=0}^{p-1}(2n+1)V_n\Ls 12^2
\e -2p^4+19p^6\mod {p^7},
\\&\sum_{n=0}^{p-1}(2n+1)V_n\Big(-\f 12\Big)^2\e p-\f 73p^4B_{p-3}\mod
{p^5},
\\&\sum_{n=0}^{p-1}(2n+1)V_n\Big(-\f 13\Big)^2
\e (-1)^{[\f p3]}p^2-10p^4U_{p-3}\mod {p^5},
\\&\sum_{n=0}^{p-1}(2n+1)V_n\Big(-\f 14\Big)^2
\e (-1)^{\f{p-1}2}p^2-10p^4E_{p-3}\mod {p^5},
\\&\sum_{n=0}^{p-1}(2n+1)V_n\Big(-\f 16\Big)^2
\e (-1)^{[\f p3]}p^2-\f{65}2p^4U_{p-3}\mod {p^5}.
\end{align*}
\vskip0.2cm
\par{\bf Conjecture 4.1} Let $p$ be an odd prime, $x\in\Bbb Z_p$, $x\not\e \pm \f 12\mod p$ and $\ls{2x-1}p=1$. If $m\in\Bbb Z^+$ and $p>4m$, then
$$\sum_{n=0}^{p-1}(2n+1)^{2m+1}v_n(x)^2\e 0\mod {p^2}.$$

\par {\bf Conjecture 4.2} Let $p>3$ be a prime. Then
\begin{align*}&
\sum_{n=0}^{p-1}(2n+1)^3V_n\Big(\f 12\Big)^2\e 6p^4-\f{143}3p^6
\mod {p^7},
\\&\sum_{n=0}^{p-1}(2n+1)^3V_n\Big(-\f 12\Big)^2\e -\f p2\mod {p^4},
\\&\sum_{n=0}^{p-1}(2n+1)^3V_n\Big(-\f 13\Big)^2\e -\f 59\Ls p3p^2\mod {p^4},
\\&\sum_{n=0}^{p-1}(2n+1)^3V_n\Big(-\f 14\Big)^2\e -\f 58(-1)^{\f{p-1}2}p^2\mod {p^4},
\\&\sum_{n=0}^{p-1}(2n+1)^3V_n\Big(-\f 16\Big)^2\e -\f {13}{18}\Ls p3p^2\mod {p^4}.
\end{align*}

\par {\bf Conjecture 4.3} Let $p>7$ be a prime. Then
\begin{align*}
&\sum_{n=0}^{p-1}(2n+1)^5V_n\Big(\f 12\Big)^2\e
\f{p}{64}\mod {p^4},
\\&\sum_{n=0}^{p-1}(2n+1)^5V_n\Big(-\f 12\Big)^2\e
\f{19}{32}p\mod {p^4},
\\&\sum_{n=0}^{p-1}(2n+1)^5V_n\Big(-\f 13\Big)^2\e \f{709}{945}\Ls p3p^2\mod {p^4},
\\&\sum_{n=0}^{p-1}(2n+1)^5V_n\Big(-\f 14\Big)^2\e \f{1841}{1920}(-1)^{\f{p-1}2}p^2\mod {p^4},
\\&\sum_{n=0}^{p-1}(2n+1)^5V_n\Big(-\f 16\Big)^2\e \f{8813}{6912}\Ls p3p^2\mod {p^4}.
\end{align*}

\section*{5. Congruences involving $g_n(x)$ and $G_n(x)$}

\par{\bf Lemma 5.1 ([4,14])} {\sl For any prime $p>3$,
\begin{align*} &G_p\e 12+64(-1)^{\f{p-1}2}p^2E_{p-3}\mod {p^3},
\\& G_p^{(3)}\e 21+243(-1)^{[\f p3]}p^2U_{p-3}\mod {p^3},\\&G_p^{(4)}\e 52+1024(-1)^{[\f p4]}p^2s_{p-3}\mod {p^3},
\\& G_p^{(6)}\e 372+8640(-1)^{\f{p-1}2}p^2E_{p-3}\mod {p^3}
\end{align*}
and
\begin{align*}&G_{p-1}\e (-1)^{\f{p-1}2}256^{p-1}+3p^2E_{p-3}\mod{p^3},
 \\&G_{p-1}^{(3)}\e (-1)^{[\f p3]}729^{p-1}+7p^2U_{p-3}\mod {p^3},
\\&G_{p-1}^{(4)}\e (-1)^{[\f p4]}4096^{p-1}+13p^2s_{p-3}\mod {p^3},
\\&G_{p-1}^{(6)}\e (-1)^{\f{p-1}2}186624^{p-1}+\f{155}9p^2E_{p-3}\mod
{p^3},\end{align*}
where $\{s_n\}$ is given by $s_0=1$ and $s_n=1-\sum_{k=0}^{n-1}\b nk2^{2n-1-2k}s_k\ (n\ge 1)$.}
 \vskip0.2cm
  \par{\bf Theorem 5.1} {\sl let $p>5$ be a prime. Then}
  \begin{align*}&
  \sum_{n=0}^{p-1}\f{G_n}{16^n}\e
  p^2\big(4(-1)^{\f{p-1}2}\cdot 16^{p-1}-3\cdot 16^{-(p-1)}\big)
  +(12-16(-1)^{\f{p-1}2})p^4E_{p-3}\mod {p^5},
  \\&\sum_{n=0}^{p-1}\f{nG_n}{16^n}
  \\&\e\f{p^2}9\big((15-3p)\cdot 16^{-(p-1)}-(16-4p)(-1)^{\f{p-1}2}
   16^{p-1}\big)
  -\f{48-80(-1)^{\f{p-1}2}}9p^4E_{p-3}
  \mod {p^5},
  \\&\sum_{n=0}^{p-1}\f{n^2G_n}{16^n}
  \\&\e -\f{p^2}{225}\big((27p^2-102p+159)\cdot 16^{-(p-1)}
  -(-1)^{\f{p-1}2}
  (36p^2-64p+112)\cdot 16^{p-1}\big)
  \\&\q+\f{16}{225}(21-53(-1)^{\f{p-1}2})p^4E_{p-3}
  \mod {p^5}.
  \end{align*}
  \par{\it Proof.} Since $\f{G_n}{16^n}=G_n(-\f 12)=g_n(\f 12)$, from Theorem 3.4 and Lemma 5.1 we deduce that
  \begin{align*}\sum_{n=0}^{p-1}\f{G_n}{16^n}
  &=\sum_{n=0}^{p-1}g_n\Ls 12=
  \f{2p^2}{\f 12-1}\Big(g_p\Ls 12-g_{p-1}\Ls 12\Big)
  =-4p^2\Big(\f{G_p}{16^p}-\f{G_{p-1}}{16^{p-1}}\Big)
  \\&\e -4p^2\Big(\f{12+64(-1)^{\f{p-1}2}p^2E_{p-3}}{16^p}
  -\f{(-1)^{\f{p-1}2}256^{p-1}+3p^2E_{p-3}}{16^{p-1}}\Big)
  \\&\e
  p^2\big(4(-1)^{\f{p-1}2}\cdot 16^{p-1}-3\cdot 16^{-(p-1)}\big)
  +(12-16(-1)^{\f{p-1}2})p^4E_{p-3}\mod {p^5}.
  \end{align*}
  Similarly, using Theorem 3.4 and Lemma 5.1 we deduce that
  \begin{align*}
  &\sum_{n=0}^{p-1}\f{nG_n}{16^n}
  \\&=\sum_{n=0}^{p-1}ng_n\Ls 12
  =\f{4p^2}9\Big((5-p)\f{G_p}{16^p}-(4-p)\f{G_{p-1}}{16^{p-1}}\Big)
  \\&\e\f{4p^2}9\Big((5-p)\f{12+64(-1)^{\f{p-1}2}p^2E_{p-3}}{16^p}
  -(4-p)\f{(-1)^{\f{p-1}2}256^{p-1}+3p^2E_{p-3}}{16^{p-1}}\Big)
  \\&\e \f{p^2}9\big((15-3p)\cdot 16^{-(p-1)}-(16-4p)(-1)^{\f{p-1}2}
   16^{p-1}\big)
  -\f{48-80(-1)^{\f{p-1}2}}9p^4E_{p-3}
  \mod {p^5}\end{align*}
  and
  \begin{align*}
  \sum_{n=0}^{p-1}\f{n^2G_n}{16^n}
  &=\sum_{n=0}^{p-1}n^2g_n\Ls 12
  =-\f{4p^2}{225}\Big((9p^2-34p+53)\f{G_p}{16^p}-(9p^2-16p+28)
  \f{G_{p-1}}{16^{p-1}}\Big)
  \\&\e \sum_{n=0}^{p-1}n^2g_n\Ls 12
  =-\f{4p^2}{225}\Big((9p^2-34p+53)\f{12+64(-1)^{\f{p-1}2}p^2E_{p-3}}
  {16^p}
  \\&\q-(9p^2-16p+28)
  \f{(-1)^{\f{p-1}2}256^{p-1}+3p^2E_{p-3}}{16^{p-1}}\Big)
  \\&  \e -\f{p^2}{225}\big((27p^2-102p+159)\cdot 16^{-(p-1)}-(-1)^{\f{p-1}2}
  (36p^2-64p+112)\cdot 16^{p-1}\big)
  \\&\q+\f{16}{225}(21-53(-1)^{\f{p-1}2})p^4E_{p-3}
  \mod {p^5}.\end{align*}
  This completes the proof.
  \vskip0.2cm
  \par In a similar way, from Theorem 3.4 and Lemma 5.1 we deduce the following congruences involving $V_n^{(r)}\ (r=3,4,6)$.
  \vskip0.2cm
  \par{\bf Theorem 5.2} {\sl } Let $p>7$ be a prime. Then
  \begin{align*}&\sum_{n=0}^{p-1}\f{G_n^{(3)}}{27^n}
    \e \f{p^2}2\big(9(-1)^{[\f p3]}\cdot 27^{p-1}-7\cdot 27^{1-p}\big)+\f{63-81(-1)^{[\f p3]}}2p^4U_{p-3}
     \mod {p^5},
  \\&\sum_{n=0}^{p-1}n\f{G_n^{(3)}}{27^n}
  \e \f{p^2}{40}\big((77-14p)\cdot 27^{1-p}-(81-18p)(-1)^{[\f p3]}\cdot 27^{p-1}\big)
  \\&\qq\qq\qq+\f{81}{40}(-7+11(-1)^{[\f p3]})p^4U_{p-3}
  \mod {p^5},
  \\&\sum_{n=0}^{p-1}n^2\f{G_n^{(3)}}{27^n}
  \e \f{p^2}{280}\big((45p^2-81p+162)(-1)^{[\f p3]}\cdot 27^{p-1}
  -(35p^2-133p+224)\cdot 27^{1-p}\big)
  \\&\qq\qq\qq+\f{81}{140}(7-16(-1)^{[\f p3]})p^4U_{p-3}
  \mod {p^5}.\end{align*}
  \par{\bf Theorem 5.3} {\sl Let $p>11$ be a prime. Then}
  \begin{align*}&\sum_{n=0}^{p-1}\f{G_n^{(4)}}{64^n}
  \e \f {p^2}3\big(16(-1)^{[\f p4]}\cdot 64^{p-1}-13\cdot 64^{1-p}\big)+\f{16}3(13-16(-1)^{[\f p4]})p^4s_{p-3}
  \mod {p^5},
  \\&\sum_{n=0}^{p-1}\f{nG_n^{(4)}}{64^n}
  \e \f {p^2}{105}\big((48p-256)(-1)^{[\f p4]}\cdot 64^{p-1}-(39p-247)\cdot 64^{1-p}\big)\\&\qq\qq\qq+\f{256}{105}(-13+19(-1)^{[\f p4]})p^4s_{p-3}
  \mod {p^5},
    \\&\sum_{n=0}^{p-1}\f{n^2G_n^{(4)}}{64^n}
  \e \f {p^2}{10395}\big(16(105p^2-192p+464)(-1)^{[\f p4]}\cdot 64^{p-1}-13(105p^2-402p+761)\cdot 64^{1-p}\big)\\&\qq\qq\qq+\f{256}{10395}(377-761(-1)^{[\f p4]})p^4s_{p-3}
  \mod {p^5}.
\end{align*}
\par{\bf Theorem 5.4} {\sl Let $p>17$ be a prime. Then}
  \begin{align*}&\sum_{n=0}^{p-1}\f{G_n^{(6)}}{432^n}
  \e \f{p^2}5\big(36(-1)^{\f{p-1}2}\cdot 432^{p-1}-31\cdot 432^{1-p}\big)+(124-144(-1)^{\f{p-1}2})p^4E_{p-3}\mod {p^5},
  \\&\sum_{n=0}^{p-1}\f{nG_n^{(6)}}{432^n}
  \e \f{p^2}{385}\Big(36(5p-36)(-1)^{\f{p-1}2}\cdot 432^{p-1}
  -31(5p-41)\cdot 432^{1-p}\Big)
  \\&\q\qq\qq+\f{144}{77}(-31+41(-1)^{\f{p-1}2})p^4E_{p-3}\mod {p^5},
    \\&\sum_{n=0}^{p-1}\f{n^2G_n^{(6)}}{432^n}
  \e \f{p^2}{85085}\Big(36(385p^2-720p+2412)(-1)^{\f{p-1}2}\cdot 432^{p-1}
  \\&\q-31(385p^2-1490p+3517)\cdot 432^{1-p}\Big)
  +\f{144}{17017}(2077-3517(-1)^{\f{p-1}2})p^4E_{p-3}
  \mod {p^5}.
  \end{align*}

  \par{\bf Lemma 5.2} {\sl Let $m,n\in\Bbb Z^+$. Then}
  \begin{align*}&\sum_{r=0}^m
  \f{(r+n+1)(r+n+2)\cdots(r+2n-1)}{(r+1)(r+2)\cdots(r+n)}
  \\&=(-1)^{n-1}\sum_{k=1}^{n-1}\b{n-1}k\b{n-1+k}k(-1)^k(H_{m+n-k}-H_{n-1-k}).
  \end{align*}
  \par{\it Proof.} Set
  $$f(x)=\f{(x+n+1)(x+n+2)\cdots(x+2n-1)}{(x+1)(x+2)\cdots(x+n)}
  =\sum_{j=1}^n\f{c_j}{x+j}.$$
  Then
  $$c_j=\lim_{x\to -j}(x+j)f(x)=\f{\prod_{1\le s\le n-1}(n+s-j)}{\prod_{1\le k\le
   n,k\not=j}(k-j)}=(-1)^{j-1}\b{n-1}{j-1}\b{2n-1-j}{n-j}.$$
Thus,
\begin{align*}
\sum_{r=0}^mf(r)&=\sum_{r=0}^m\sum_{j=1}^n\f{c_j}{r+j}
=\sum_{j=1}^nc_j\sum_{r=0}^m\f 1{r+j}
\\&=\sum_{j=1}^n(-1)^{j-1}\b{n-1}{j-1}\b{2n-1-j}{n-j}
(H_{m+j}-H_{j-1})
\\&=\sum_{k=1}^n(-1)^{n-k-1}\b{n-1}{n-k-1}\b{2n-1-(n-k)}k
(H_{m+n-k}-H_{n-k-1}),
\end{align*}
which yields the result.
\vskip0.2cm
  \par{\bf Lemma 5.3} {\sl Let $p$ be an odd prime and $x\in\Bbb Z_p$ with $\ls{2x-1}p=1$. Then $g_{p-1}(x),g'_{p-1}(x),g_p(x),$ $pg_p'(x)\in\Bbb Z_p$.
Moreover, }
\begin{align*}&g_p(x)\e
 \Big(1+\f{x_1(x_1+8x_0-4)}{16(2x-1)}\Big)(1+2pH_{x_0-1})
+p\f{x_1^2(x_1+12x_0-6)}{32(2x-1)^2}
\mod {p^2},
\\&pg'_p(x)\e
\f{x_1+4x_0-2}{4(2x-1)}
\Big(1+p\Big(\f {x_1(x_1+8x_0-4)}{4(2x-1)(x_1+4x_0-2)}+2H_{x_0-1}\Big)\Big)
\mod {p^2},
\\&g_{p-1}(x)\e
(-1)^{x_0-1}(1-2pH_{x_0-1})
\mod {p^2},
 \\&g_{p-1}'(x)\e 0\mod p.
\end{align*}
where $x_0$ and $x_1$ are given by $x_0\in\{1,2,\ldots,\f{p-1}2\}$, $2x-1\e (2x_0-1)^2\mod p$ and $2x-1=(2x_0-1)^2+px_1$.
\vskip0.2cm
\par{\it Proof.}
  For $k=1,2,\ldots,p-1$ we see that
\begin{align*}
&\f{\prod_{r=1}^k(2x-1-(2r-1)^2)}{4^k
\cdot k!^2}
\\&\e \f{\prod_{r=1}^k((2x_0-1)^2-(2r-1)^2)}{4^k
\cdot k!^2}=\f{\prod_{r=1}^k(x_0-r)(x_0+r-1)}{k!^2}
\\&=\begin{cases}
\b{x_0-1}k\b{x_0-1+k}k\mod p&\t{if $k\le x_0-1$,}
\\0\mod p&\t{if $k\ge x_0$}.\end{cases}
\end{align*}
By (4.1),
$$\sum_{k=1}^n\b nk\b{n+k}k\f{(-1)^k}k=-2H_n.$$
Since $\b pk=\f pk\b{p-1}{k-1}\e -(-1)^k\f pk\mod {p^2}$ for $1\le k\le p-1$,
from the above we deduce that
\begin{align*}
&\sum_{k=1}^{p-1}\b pk\f{\prod_{r=1}^k(2x-1-(2r-1)^2}{4^k\cdot k!^2}
\\&\e -p\sum_{k=1}^{x_0-1}\b{x_0-1}k\b{x_0-1+k}k\f{(-1)^k}k
=2pH_{x_0-1}\mod {p^2}.
\end{align*}
Thus,
\begin{align*}
g_p(x)&=1+\sum_{k=1}^{p}\b{p}k\f{\prod_{r=1}^k(2x-1
-(2r-1)^2)}{4^k
\cdot k!^2}
\\&\e 1+2pH_{x_0-1}+\f{\prod_{r=1}^p(2x-1-(2r-1)^2)}{4^p
\cdot p!^2}
\mod {p^2}.
\end{align*}
By the congruence for
$\prod_{r=1}^p(2x-1-(2r-1)^2)$
in the proof of Lemma 4.1,
\begin{align*}\prod_{r=1}^p(2x-1-(2r-1)^2)
&\e p^2x_1(x_1+4(2x_0-1)-4p)\cdot \f{4^{p-1}\cdot (p-1)!^2}{4(2x-1)}
\\&\q\times\Big(1+p\Big(\f 1{2x_0-1}+2H_{x_0-1}+\f {x_1}{2(2x-1)}\Big)\Big)\mod {p^4}.
\end{align*}
Thus,
\begin{align*}
g_p(x)&\e 1
+2pH_{x_0-1}+\f{\prod_{r=1}^p(2x-1-(2r-1)^2)}{4^p
\cdot p!^2}
\\&\e 1+2pH_{x_0-1}+
x_1(x_1+4(2x_0-1)-4p)\cdot \f{1}{16(2x-1)}
\\&\q\times\Big(1+p\Big(\f 1{2x_0-1}+2H_{x_0-1}+\f {x_1}{2(2x-1)}\Big)\Big)
\\&\e \Big(1+\f{x_1(x_1+8x_0-4)}{16(2x-1)}\Big)(1+2pH_{x_0-1})
+p\f{x_1^2(x_1+12x_0-6)}{32(2x-1)^2}
\mod {p^2}.
\end{align*}
Since
$$g_p(t)=1+\sum_{k=1}^p\b pk\f{\prod_{r=1}^k(2t-1-(2r-1)^2)}{4^k
\cdot k!^2},$$
  we see that
  \begin{align*}g'_p(x)&=\sum_{k=1}^p\b pk\f{\sum_{s=1}^k2\prod_{1\le r\le k,r\not=s}(2x-1-(2r-1)^2)}{4^k\cdot k!^2}
\\&=\f 12\sum_{k=1}^{p-1}\b pk\f{\sum_{s=1}^k\prod_{1\le r\le k,r\not=s}(px_1+(2x_0-1)^2-(2r-1)^2)}{4^{k-1}\cdot k!^2}
\\&\q+\f 12\cdot \f{\sum_{s=1}^p\prod_{1\le r\le p,r\not=s}(px_1+(2x_0-1)^2-(2r-1)^2)}{4^{p-1}\cdot p!^2}.
\end{align*}
By the formula for $\f 12f'_p(x)$ in the proof of Lemma 4.1,
\begin{align*}
&\sum_{s=1}^p\prod_{1\le r\le p,r\not=s}(px_1+(2x_0-1)^2-(2r-1)^2)
\\&\e  p\cdot \f{4^{p-1}\cdot (p-1)!^2}{4(2x-1)}
\Big(1+p\Big(\f 1{2x_0-1}+2H_{x_0-1}+\f {x_1}{2(2x-1)}\Big)\Big)
\\&\q\times \Big(2x_1+4(2x_0-1)-4p
-\f{px_1(x_1+4(2x_0-1))}{2(2x-1)}\Big)
\mod {p^3}.
\end{align*}
Note that $p\mid \b pk$ for $k=1,2,\ldots,p-1$. From the above we deduce that
\begin{align*}2pg'_p(x)&
\e p\f{\sum_{s=1}^p\prod_{1\le r\le p,r\not=s}(px_1+(2x_0-1)^2-(2r-1)^2)}{4^{p-1}\cdot p!^2}
\\&\e \f 1{4(2x-1)}
\Big(1+p\Big(\f 1{2x_0-1}+2H_{x_0-1}+\f {x_1}{2(2x-1)}\Big)\Big)
\\&\q\times \Big(2x_1+8x_0-4-4p
-\f{px_1(x_1+8x_0-4)}{2(2x-1)}\Big)
\\&\e \f{x_1+4x_0-2}{2(2x-1)}
\Big(1+p\Big(\f {x_1(x_1+8x_0-4)}{4(2x-1)(x_1+4x_0-2)}+2H_{x_0-1}\Big)\Big)
\mod {p^2}.
\end{align*}
\par Now we deduce the congruences for $g'_{p-1}(x)$ modulo $p$ and $g_{p-1}(x)$ modulo $p^2$. Since
$$g_{p-1}(x)=1+\sum_{k=1}^{p-1}\b{p-1}k
\f{\prod_{r=1}^k(2x-1-(2r-1)^2)}{4^k
\cdot k!^2},$$
$(2x_0-1)^2\e (2r-1)^2\mod p$ for $r=x_0,p+1-x_0$ and $\b{p-1}k=\f{(p-1)(p-2)\cdots(p-k)}{k!}\e (-1)^k(1-pH_k)
\mod{p^2}$ for $1\le k\le p-1$, applying Vandermonde's identity we derive
\begin{align*}
g_{p-1}(x)&=1+\sum_{k=1}^{p-1}\b{p-1}k
\f{\prod_{r=1}^k(px_1+(2x_0-1)^2
-(2r-1)^2)}{4^k \cdot k!^2}
\\&\e 1+\sum_{k=1}^{x_0-1}\b{p-1}k
\f{\prod_{r=1}^k((2x_0-1)^2-(2r-1)^2)\cdot(1+\sum_{s=1}^k
\f{px_1}{(2x_0-1)^2-(2s-1)^2})}{4^k\cdot k!^2}
\\&\q+\sum_{k=x_0}^{p-x_0}\b {p-1}k
\f{px_1\prod_{1\le r\le k,r\not=x_0}((2x_0-1)^2-(2r-1)^2)}{4^k\cdot k!^2}
\\&\e 1+\sum_{k=1}^{x_0-1}(1-pH_k)(-1)^k\f{\prod_{r=1}^k
(x_0-r)(x_0+r-1)}{k!^2}\Big(1+\sum_{s=1}^k
\f{px_1}{4(x_0-s)(x_0+s-1)}\Big)
\\&\q+\sum_{k=x_0}^{p-x_0}(-1)^k\f{px_1}4
\cdot
\f{\prod_{1\le r\le k,r\not=x_0}(x_0-r)(x_0+r-1)}{k!^2}
\\&\e 1+\sum_{k=1}^{x_0-1}(-1)^k\b{x_0-1}k\b{x_0-1+k}k
-p\sum_{k=1}^{x_0-1}(-1)^kH_k\b{x_0-1}k\b{x_0-1+k}k
\\&\q+\f{px_1}4\sum_{k=1}^{x_0-1}(-1)^k\b{x_0-1}k\b{x_0-1+k}k
\cdot \f 1{2x_0-1}\sum_{s=1}^k\Big(\f 1{x_0-s}+\f 1{x_0+s-1}\Big)
\\&\q+\f {px_1}4\sum_{k=x_0}^{p-x_0}(-1)^{x_0}
\f{(x_0-1+k)!(k-x_0)!}{k!^2\cdot (2x_0-1)}
\\&=\sum_{k=0}^{x_0-1}\b{x_0-1}{x_0-1-k}\b{-x_0}k
-p\sum_{k=1}^{x_0-1}\b{x_0-1}k\b{x_0-1+k}k(-1)^kH_k
\\&\q+\f{px_1}{4(2x_0-1)}
\Big(\sum_{k=1}^{x_0-1}\b{x_0-1}k\b{x_0-1+k}k(-1)^k
(H_{x_0-1+k}-H_{x_0-1-k})
\\&\q+(-1)^{x_0}\sum_{k=x_0}^{p-x_0}
\f{(k+1)(k+2)\cdots(k+x_0-1)}{k(k-1)\cdots(k-x_0+1)}\Big)
\\&=\b{-1}{x_0-1}-p\sum_{k=1}^{x_0-1}\b{x_0-1}k\b{x_0-1+k}k(-1)^kH_k
\\&\q+\f{px_1}{4(2x_0-1)}
\Big(\sum_{k=1}^{x_0-1}\b{x_0-1}k\b{x_0-1+k}k(-1)^k
(H_{x_0-1+k}-H_{x_0-1-k})
\\&\q+(-1)^{x_0}\sum_{r=0}^{p-2x_0}
\f{(r+x_0+1)(r+x_0+2)\cdots(r+2x_0-1)}{(r+1)
(r+2)\cdots(r+x_0)}\Big)
\mod {p^2}.
\end{align*}
It is well known (see for example [15, p.242]) that
$\sum_{k=0}^n\b nk(-1)^kH_k=-\f 1n$ for $n\ge 1$. Thus,
applying the formula in [15, p.259] and (4.1) gives
$$\sum_{k=1}^n\b nk\b{n+k}k(-1)^kH_{k}
=-(-1)^n\sum_{k=1}^n\b nk\b{n+k}k\f{(-1)^k}k=2(-1)^nH_n.\eqno{(5.1)}$$
From Lemma 5.2 and the fact that $H_{p-1-n}
=\sum_{k=1}^{p-1}\f 1k-\sum_{k=1}^n\f 1{p-k}\e H_n\mod p$ we see that
\begin{align*}
&\sum_{r=0}^{p-2x_0}
\f{(r+x_0+1)(r+x_0+2)\cdots(r+2x_0-1)}{(r+1)
(r+2)\cdots(r+x_0)}
\\&=(-1)^{x_0-1}\sum_{k=1}^{x_0-1}\b{x_0-1}k\b{x_0-1+k}k(-1)^k
(H_{p-2x_0+x_0-k}-H_{x_0-1-k})
\\&\e (-1)^{x_0-1}\sum_{k=1}^{x_0-1}\b{x_0-1}k\b{x_0-1+k}k(-1)^k
(H_{x_0-1+k}-H_{x_0-1-k})\mod p.\tag{5.2}
\end{align*}
Thus,
\begin{align*}g_{p-1}(x)&\e \b{-1}{x_0-1}-p\sum_{k=1}^{x_0-1}\b{x_0-1}k\b{x_0-1+k}k(-1)^kH_k
\\&\e (-1)^{x_0-1}-2p(-1)^{x_0-1}H_{x_0-1}\mod {p^2}.
\end{align*}
On the other hand,
\begin{align*} g'_{p-1}(x)&=\sum_{k=1}^{p-1}\b{p-1}k
\f{\sum_{s=1}^k 2\prod_{1\le r\le k,r\not=s}(2x-1-(2r-1)^2)}{4^k\cdot k!^2}
\\&\e \sum_{k=1}^{p-1}\f {(-1)^k}{2\cdot 4^{k-1}\cdot k!^2}
\sum_{s=1}^k\prod_{1\le r\le k,r\not=s}((2x_0-1)^2-(2r-1)^2)
 \\&\e \f 12\sum_{k=1}^{x_0-1}\f {(-1)^k}{4^{k-1}\cdot k!^2}
\sum_{s=1}^k\prod_{1\le r\le k,r\not=s}((2x_0-1)^2-(2r-1)^2)
\\&\q+\f 12\sum_{k=x_0}^{p-x_0}\f {(-1)^k}{4^{k-1}\cdot k!^2}
\prod_{1\le r\le k,r\not=x_0}((2x_0-1)^2-(2r-1)^2)
\\&\e \f 12\sum_{k=1}^{x_0-1}\f {(-1)^k}{ k!^2}
\prod_{1\le r\le k}(x_0-r)(x_0+r-1)\cdot\Big(\sum_{s=1}^k
\f 1{(x_0-s)(x_0+s-1)}\Big)
\\&\q+\f 12\sum_{k=x_0}^{p-x_0}\f {(-1)^k}{k!^2}
\prod_{1\le r\le k,r\not=x_0}(x_0-r)(x_0+r-1)
\\&=\f 1{2(2x_0-1)}\Big(\sum_{k=1}^{x_0-1}\b{x_0-1}k\b{x_0-1+k}k(-1)^k
\Big(\sum_{s=1}^k\f 1{x_0-s}+\sum_{s=1}^k\f 1{x_0+s-1}\Big)
\\&\q+\sum_{k=x_0}^{p-x_0}(-1)^{x_0}\f{(x_0-1+k)!(k-x_0)!}{k!^2}
\Big)
\\&=\f 1{2(2x_0-1)}\Big(\sum_{k=1}^{x_0-1}\b{x_0-1}k\b{x_0-1+k}k(-1)^k
(H_{x_0-1+k}-H_{x_0-1-k})
\\&\q+(-1)^{x_0}\sum_{r=0}^{p-2x_0}
\f{(r+x_0+1)(r+x_0+2)\cdots(r+2x_0-1)}{(r+1)
(r+2)\cdots(r+x_0)}\Big)
\mod p.\end{align*}
This together with (5.2) yields $g'_{p-1}(x)\e 0\mod p$.
The proof is now complete.
\vskip0.2cm
\par{\bf Remark 5.1} In order to prove Lemma 5.3, we find the following identities:
\begin{align*}
&\sum_{k=0}^n\b nk\b{n+k}k(-1)^kH_{n+k}=2(-1)^nH_n,\tag{5.3}
\\&\sum_{k=0}^n\b nk\b{n+k}k(-1)^kH_{n-k}=(-1)^{n-1}\sum_{k=1}^n\b nk^2\f 1k.
\tag {5.4}
\end{align*}
 \par{\bf Theorem 5.5} {\sl Let $p$ be an odd prime and $x\in\Bbb Z_p$ with $\ls{2x-1}p=1$. Suppose that $x_0$ and $x_1$ are given by $x_0\in\{1,2,\ldots,\f{p-1}2\}$, $2x-1\e (2x_0-1)^2\mod p$ and $2x-1=(2x_0-1)^2+px_1$.
 \par $(\t{\rm i})$ For $x\not\e 1\mod p$,
 \begin{align*}\sum_{n=0}^{p-1}g_n(x)
 &\e  \f{2p^2}{x-1}\Big(1-(-1)^{x_0-1}+\f{x_1(x_1+8x_0-4)}
 {16(2x-1)}\Big)
 \\&\q+\f{4p^3}{x-1}\Big(1+(-1)^{x_0-1}+\f{x_1(x_1+8x_0-4)}
 {16(2x-1)}\Big)H_{x_0-1} \mod {p^4}.
 \end{align*}
 \par $(\t{\rm ii})$ For $x\not\e 1,5\mod p$,
 \begin{align*}
 &\sum_{n=0}^{p-1}ng_n(x)
 \\&\e -\f{p^2}{(x-1)(x-5)}\Big(2x-6+4(-1)^{x_0-1}+
 \f{(x-3)x_1(x_1+8x_0-4)}
 {8(2x-1)}\Big)
 \\&\q+\f{2p^3}{(x-1)(x-5)}
 \Big((x-1)\Big(1-(-1)^{x_0-1}+\f{x_1(x_1+8x_0-4)}{16(2x-1)}\Big)
 -\f{(x-3)x_1^2(x_1+12x_0-6)}{32(2x-1)^2}
 \\&\q+\Big(4(-1)^{x_0-1}-2x+6-\f{(x-3)x_1(x_1+8x_0-4)}{8(2x-1)}
 \Big) H_{x_0-1}\Big)\mod {p^4}.
  \end{align*}
 \par $(\t{\rm iii})$ For $x\not\e 1,5,13\mod p$,}
 \begin{align*}&\sum_{n=0}^{p-1}n^2g_n(x)
 \\&\e \f{p^2}{(x-1)(x-5)(x-13)}
 \Big(2(x^2-12x+19)-4(-1)^{x_0-1}(x+3)
 \\&\q+(x^2-12x+19)\f{x_1(x_1+8x_0-4)}{8(2x-1)}\Big)
 \\&\q+\f{2p^3}{(x-1)(x-5)(x-13)}\Big(-2(x-1)(x-9)
 \Big(1+\f{x_1(x_1+8x_0-4)}{16(2x-1)}\Big)
 \\&\q+(x^2-12x+19)\f{x_1^2(x_1+12x_0-6)}{32(2x-1)^2}-8(-1)^{x_0-1}(x-1)
 \\&\q+\Big(2(x^2-12x+19)\Big(1+\f{x_1(x_1+8x_0-4)}{16(2x-1)}\Big)+4(-1)^{x_0-1}
 (x+3)\Big)H_{x_0-1}\Big)\mod {p^4}.
 \end{align*}
 \vskip0.2cm
\par{\it Proof.} This is immediate from Theorem 3.4 and Lemma 5.3.
\vskip0.2cm
\par{\bf Theorem 5.6} {\sl Let $p$ be an odd prime, $x\in\Bbb Z_p$ and $\ls{2x-1}p=1$. Then
\begin{align*}\sum_{n=0}^{p-1}g_n(x)^2&\e
(-1)^{x_0-1}\f{x_1+4x_0-2}{2(2x-1)}p
+(-1)^{x_0-1}\f{x_1(x_1+8x_0-4)}{8(2x-1)^2}p^2
\mod {p^3},\end{align*}
where $x_0$ and $x_1$ are given by $x_0\in\{1,2,\ldots,\f{p-1}2\}$, $2x-1\e (2x_0-1)^2\mod p$ and $2x-1=(2x_0-1)^2+px_1$.}
\vskip0.2cm
\par{\it Proof.} From (3.11) and Lemma 5.3,
\begin{align*}\sum_{n=0}^{p-1}g_n(x)^2
&=2p\cdot g_{p-1}(x)\cdot pg'_p(x)-2p^2g_p(x)g'_{p-1}(x)
\\&\e2p\cdot (-1)^{x_0-1}(1-2pH_{x_0-1})
\\&\q\times \f{x_1+4x_0-2}{4(2x-1)}
\Big(1+p\Big(\f {x_1(x_1+8x_0-4)}{4(2x-1)(x_1+4x_0-2)}+2H_{x_0-1}\Big)\Big)
\mod {p^3}, \end{align*}
which yields the result.
\vskip0.2cm
\par{\bf Corollary 5.1} Let $p$ be an odd prime, $x\in\Bbb Z_p$ and $x\not\e -\f 12\mod p$. Then
$$\sum_{n=0}^{p-1}G_n(x)^2\e (-1)^{\xp}\f{1+2x'}{1+2x}p\mod
{p^3},\eqno{(5.5)} $$ where $\xp\in\{0,1,\ldots,p-1\}$ is given by $x\e \xp\mod p$ and $x'=(x-\xp)/p$.
\vskip0.2cm
\par{\it Proof.} Observe that $G_n(x)=g_n(2x^2+2x+1)$ and
$$2(2x^2+2x+1)-1=(2x+1)^2\e (2(\xp+1)-1)^2
\e (2(p-\xp)-1)^2\mod p.$$
Set
$$ x_0=\begin{cases}\xp+1=x+1-px'&\t{if $\xp<\f{p-1}2$,}
\\p-\xp=(x'+1)p-x&\t{if $\xp\ge \f{p-1}2$}
\end{cases}$$
and
$$ x_1=\f{(2x+1)^2-(2x_0-1)^2}p=\begin{cases}
4x'(2x+1-px')&\t{if $\xp<\f{p-1}2$,}
\\ 4(x'+1)(2x+1-p(x'+1))&\t{if $\xp\ge \f{p-1}2$.}
\end{cases}$$
Then $x_0\in\{1,2,\ldots,\f{p-1}2\}$ and
$2(2x^2+2x+1)-1=(2x_0-1)^2+px_1$. From Theorem 5.6 we have
 \begin{align*}\sum_{n=0}^{p-1}G_n(x)^2&=
 \sum_{n=0}^{p-1}g_n(2x^2+2x+1)^2
 \\&\e
(-1)^{x_0-1}\f{x_1+4x_0-2}{2(2x+1)^2}p
+(-1)^{x_0-1}\f{x_1(x_1+8x_0-4)}{8(2x+1)^4}p^2
\mod {p^3}.\end{align*}
It is easy to see that
$$x_1+4x_0-2=2(2x+1)(2x'+1)-4px'(x'+1)$$
and
$$x_1+8x_0-4\e 16(2x+1)^2x'(x'+1)\mod p.$$
Thus,
\begin{align*}
\sum_{n=0}^{p-1}G_n(x)^2
&\e (-1)^{\xp}
\f{2(2x+1)(2x'+1)-4px'(x'+1)}{2(2x+1)^2}p
+(-1)^{\xp}\f{16(2x+1)^2x'(x'+1)}{8(2x+1)^4}p^2
\\&=(-1)^{\xp}\f{2x'+1}{2x+1}p\mod {p^3}.
\end{align*}
This proves the corollary.
\vskip0.2cm
\par We note that Corollary 5.1 was first conjectured by Z.W. Sun[17], and he proved the congruence modulo $p^2$.
\vskip0.2cm
\par In [14] the author obtained the congruence for $G_{\f{p-1}2}$ modulo $p^2$ for any prime $p>3$. Calculation with Maple suggests the following conjecture.
\vskip0.2cm
\par{\bf Conjecture 5.1} For any prime $p>3$,
$$G_{\f{p-1}2}\e\begin{cases} 2^{p-1}\b{(p-1)/2}{(p-1)/4}^2\big(1-\f{p^2}6E_{p-3}\big)\pmod{p^3}
&\t{if $p\e 1\pmod 4$,}
\\\f{p^2}{3\b{(p-1)/2}{(p-3)/4}^2}\pmod {p^3}&\t{if $p\e 3\pmod 4$.}
\end{cases}$$
\par  Let $x$ be a real number and $n\in\Bbb Z^+$. In [14] and [16] the author conjectured that
$G_n(x)^2<G_{n+1}(x)G_{n-1}(x)$ and
 $V_n(x)^2<V_{n+1}(x)V_{n-1}(x)$ for $x\in (-1,0)$, and that $G_n(x)^2>G_{n+1}(x)G_{n-1}(x)$ and $V_n(x)^2>V_{n+1}(x)V_{n-1}(x)$ for $x\notin [-1,0]$.
This has been confirmed by Mao and Xiao in [9]. Now we pose new challenging conjecture for $g_n(x)$ and $v_n(x)$.
\vskip0.2cm
\par{\bf Conjecture 5.2} {\sl Let $n$ be a positive integer. Then $g_n(x)^2>g_{n-1}(x)g_{n+1}(x)$ for $x\le -1$, and
$v_n(x)^2>v_{n+1}(x)v_{n-1}(x)$ for $x\le -\f 18$.}


\begin{thebibliography}{99}
\bibitem [1] {} G. Andrews, R. Askey and R. Roy, {\it Special Functions}, Cambridge Univ. Press, New York, 1999.

    \bibitem [2] {} P. Borwein and T. Erd\'elyi, {\it
    Polynomials and Polynomial Inequalities}, Springer, New York, 1995.

\bibitem [3] {} Y.-T. Li, X.-S. Wang and R. Wong, {\it Asymptotics of the Wilson polynomials}, Anal. Appl. {\bf 18}(2020), 237-270.

    \bibitem [4] {} J.-C. Liu and H.-X. Ni, {\it On two supercongruences involving
Almkvist-Zudilin sequences}, Czechoslovak Math. J. {\bf 71}(2021), 1211-1219.

\bibitem [5] {} W. Magnus, F. Oberhettinger and R. P. Soni,
{\it Formulas and Theorems for the Special Functions of Mathematical Physics}, 3rd edn.
Springer, New York, 1966, p.206.

\bibitem [6] {}  G.-S. Mao, {\it On some conjectural congruences involving Ap\'ery-like numbers $V_n$}, preprint (2021), ResearchGate: 10.13140/RG.2.2.12308.42880.

\bibitem [7] {} G.-S. Mao, {\it Proof of some conjectural congruences involving Ap\'ery-like sequences}, J. Difference Equ. Appl. {\bf 31}(2025), 570-587.

         \bibitem [8] {} G.-S. Mao and J.-J. Yang,
       {\it Supercongruences involving Ap\'ery-like sequences $\{m^nG_n(x)\}$
and $\{m^nV_n(x)\}$}, Period. Math. Hung. 87(2023), 303-314.

\bibitem [9] {} J. Mao and Q. Xiao, {\it Second order determinants for three-term recurrence polynomial sequences}, Filomat {\bf 40}(2026), 3169-3178.

\bibitem [10] {} N.J.A. Sloane, {\it The On-Line Encyclopedia
of Integer Sequences}, http://oeis.org/.

\bibitem [11] {} Z.H. Sun, {\it Congruences involving Bernoulli and Euler numbers}, J. Number Theory 128(2008), 280-312.

\bibitem [12] {} Z.H. Sun, {\it
Super congruences concerning Bernoulli polynomials}, Int. J. Number
Theory {\bf 11}(2015), 2393-2404.

\bibitem [13] {} Z.H. Sun, {\it Congruences involving binomial
coefficients and Ap\'ery-like numbers}, Publ. Math. Debrecen {\bf
96}(2020), 315-346.

\bibitem [14] {} Z.H. Sun, {\it
Congruences for certain families of Ap\'ery-like sequences}, Czech.
Math. J. {\bf 72}(2022), 875-912.

\bibitem [15] {}  Z.H. Sun, {\it Binomial Coefficients, Recurrence Sequences and Congruences} (Chinese), Science Press, Beijing, 2025.

\bibitem [16] {} Z. H. Sun, {\it Congruences for a type of Ap\'ery-like numbers}, Chin. Ann. Math., Ser. B, to appear.

    \bibitem [17] {} Z.W. Sun, {\it
Supercongruences involving dual sequences}, Finite Fields Appl. {\bf
46}(2017), 179-216.


\bibitem [18] {} S. Yang and J.-C. Liu, {\it On the divisibility of sums involving Ap\'ery-like polynomials},
    Bull. Aust. Math. Soc. {\bf 106}(2022), 203-208.


\bibitem [19] {} C. Wang, {\it On two conjectural supercongruences of Z.-W. Sun}, Ramanujan J. {\bf 56}(2021), 1111-1121.

\end{thebibliography}
\end{document}